\documentclass[onefignum,onetabnum]{siamart190516}

\usepackage[utf8]{inputenc}
\usepackage[T1]{fontenc}
\usepackage{color}
\usepackage{csquotes} 
\usepackage[caption=false]{subfig}
\usepackage{microtype} 
\usepackage[normalem]{ulem}
\usepackage{soul}
\usepackage{siunitx}
\usepackage{dsfont} 
\usepackage{mathtools}
\usepackage{xpatch} 
\xpretocmd{\eqref}{Eq.~}{}{} 
\usepackage{cases} 
\usepackage{algorithm}
\usepackage{algpseudocode}
\usepackage{enumitem} 
\usepackage{physics}
\usepackage{xspace}
\usepackage{makecell}
\usepackage{booktabs} 
\usepackage{tabularx}
\usepackage{pgfplots}
\usepackage{pgfplotstable}

\usepackage{amsfonts}

\newsiamremark{remark}{Remark}
\newsiamremark{hypothesis}{Hypothesis}
\crefname{hypothesis}{Hypothesis}{Hypotheses}
\newsiamthm{claim}{Claim}

\pgfplotsset{
	compat=1.13,
}
\usetikzlibrary{external,pgfplots.groupplots}
\tikzsetexternalprefix{figures_ext/}

\newcommand{\tool}[1]{#1}

\newcommand*{\wina}{Winkelmann et al.\xspace}

\newcommand*{\BFGS}{\tool{BFGS}}
\newcommand*{\FEAST}{\tool{FEAST}}

\newcommand*{\LBFGSB}{\tool{L-BFGS-B}}
\newcommand*{\SLISEtext}{SLiSe}
\newcommand*{\SLISE}{\tool{\SLISEtext{}}}
\newcommand*{\WCRtext}{WCR}
\newcommand*{\WCR}{\tool{\WCRtext{}}}

\newcommand*{\ESlise}{WiSe}

\newcommand*{\C}{\mathbb{C}}
\newcommand*{\R}{\mathbb{R}}
\newcommand*{\N}{\mathbb{N}}

\newcommand*{\HR}{\mathbb{H}^{+R}}
\newcommand*{\eps}{\varepsilon}
\newcommand*{\iu}{{i\mkern1mu}}
\newcommand*{\wgtf}{\omega}
\newcommand*{\wgtfsearch}[1]{\wgtf_{#1}}
\newcommand*{\wgtfgamma}{\wgtf_{\gamma\text{-\SLISE{}}}}
\newcommand*{\wgtfs}{V_s}
\newcommand*{\ind}{\mathds{1}_{(-1,1)}}
\newcommand*\indi[1]{\mathds{1}_{(#1)}}

\DeclareMathOperator{\sign}{sign}

\DeclareMathOperator{\am}{argmin}

\newcommand*\argmin[2]{\underset{#1}{\am} \ {#2}}
\newcommand{\set}[1]{\{\,#1\,\}}
\newcommand*{\abbs}[1]{\left\lvert#1\right\rvert}

\newcommand*\colvecalt[2]{\begin{pmatrix}#1 \\ #2\end{pmatrix}}
\newcommand*\conj[1]{\overline{#1}}

\newcommand*\ratcp[2]{#1_{#2}}
\newcommand{\expnumber}[2]{{#1}\mathrm{e}{#2}}

\mathchardef\ordinarycolon\mathcode`\:
\mathcode`\:=\string"8000
\begingroup \catcode`\:=\active
\gdef:{\mathrel{\mathop\ordinarycolon}}
\endgroup

\setdescription{style=sameline,leftmargin=0.6in} 

\algdef{SE}[DOWHILE]{Do}{doWhile}{\algorithmicdo}[1]{\algorithmicwhile\ #1}%

\makeatletter
\newcommand{\algmargin}{\the\ALG@thistlm}
\makeatother
\newlength{\whilewidth}
\algdef{SE}[parWHILE]{parWhile}{EndparWhile}[1]
{\parbox[t]{\dimexpr\linewidth-\algmargin}{%
		\hangindent\whilewidth\strut\algorithmicwhile\ #1\ 
		\algorithmicdo\strut}}{\algorithmicend\ \algorithmicwhile}%
\algnewcommand{\parState}[1]{\State%
	\parbox[t]{\dimexpr\linewidth-\algmargin}{\strut #1\strut}}

\definecolor{SLISEC}{RGB}{77,175,74} 
\definecolor{GAUSSC}{RGB}{55,126,184} 
\definecolor{ZOLOC}{RGB}{228,26,28} 
\definecolor{SLISEC2}{RGB}{152,78,163} 
\definecolor{SLISEC3}{RGB}{255,127,0} 
\definecolor{XAXIS}{RGB}{100,100,100}

\newcolumntype{D}{>{\centering\arraybackslash}p{1.6cm}} 

\newcommand\cmt[2]{
	\noindent\ignorespaces \textcolor{blue}\bgroup
	\textbf{#1} -- Edo: %
	\textit{#2}\egroup
}

\definecolor{brown}{rgb}{.5,.15,.15}
\usepackage{setspace}

\newcommand\kmt[1]{
	\textcolor{blue}\bgroup KK:
	\em \textbf{#1}
	\egroup
}

\headers{Rational spectral filters with optimal convergence rate}{K. Kollnig, P. Bientinesi, and E. Di Napoli}

\title{Rational spectral filters with optimal convergence rate
}

\author{Konrad Kollnig\thanks{RWTH Aachen University
  (\email{konrad.kollnig@rwth-aachen.de}).}
\and Paolo Bientinesi\thanks{Ume\aa\ University
(\email{pauldj@cs.umu.se}).}
\and Edoardo {Di Napoli}\thanks{Jülich Supercomputing Centre, Forschungszentrum Jülich
(\email{e.di.napoli@fz-juelich.de}).}}

\usepackage{amsopn}

\newcommand*\ratslise[4]{\ratcp{#1}{#2, #3}(#4)}

\ifpdf
\hypersetup{
  pdftitle={Rational spectral filters with optimal convergence rate},
  pdfauthor={K. Kollnig, P. Bientinesi, and E. Di Napoli}
}
\fi

\begin{document}

\maketitle

\begin{abstract}
  In recent years, contour-based eigensolvers have emerged as a
  standard approach for the solution of large and sparse eigenvalue
  problems. Building upon recent performance improvements through
  non-linear least square optimization of so-called rational filters,
  we introduce a systematic method to design these filters by
  minimizing the worst-case convergence ratio and eliminate the
  parametric dependence on weight functions. Further, we provide an
  efficient way to deal with the box-constraints which play a central
  role for the use of iterative linear solvers in contour-based
  eigensolvers. Indeed, these parameter-free filters consistently
  minimize the number of iterations and the number of FLOPs to reach
  convergence in the eigensolver. As a byproduct, our rational filters
  allow for a simple solution to load balancing when the solution of
  an interior eigenproblem is approached by the slicing of the sought
  after spectral interval.
\end{abstract}


\begin{keywords}
  Hermitian Eigenvalue Problem, Rational Filters, Contour-based Eigensolver, 
  FEAST, Worst-case Convergence Rate, Load Balancing, Non-linear Least
  Squares, BFGS, Nelder-Mead.
\end{keywords}

\begin{AMS}
	65F15, 41A20, 65Y05
\end{AMS}

\section{Introduction}
\label{sec:introduction}
For the Hermitian eigenproblem $Ax = \lambda x$ with $\lambda \in [a,b] \subset \R$, 
the last decade has seen the emergence of a new class of eigensolvers
based on spectral projectors. Such eigensolvers are typically expressed as
integrals of the spectral resolvent $(A - zI)^{-1}$ over a contour in the complex
plane that encloses the interval $[a,b]$ 
\cite{sakurai,sakurai_cirr,feast,ifeast,slicing_library}.
Numerical quadrature transforms the contour integral into a
matrix-valued rational function with complex coefficients $\beta_i$
and poles $z_i$.  In this form, the problem of finding an efficient
spectral projector is mapped to that of finding a rational function---often
referred to as rational filter---that approximates the indicator function
\begin{equation}
 \label{eqn:indicator}
 \indi{a,b}(x) = \begin{cases} 1, & \text{if } x \in [a,b], \\
  0, & \text{otherwise.}\end{cases}
\end{equation}
This is a discontinuous function, often termed the ``ideal filter'',
because it exactly maps the desired eigenvalues in the interval to $1$ and the
rest of the spectrum to $0$.

The algorithmic structure of eigensolvers based on rational filters
has the advantage of lending itself to parallel implementations with
multiple levels of nested parallelism~\cite{zparsesolver,feastsolver}. On the other hand, several
factors make load balancing for these parallel eigensolvers a potential
nightmare~\cite{goeddeke,guettel}. Among them the design of the filter is an important element
that influences the convergence of the eigensolver with direct
consequences on the load balancing of any parallel implementation
based on slicing $[a, b]$ in subintervals.
In this paper, we focus on the design of filters with the aim of
resolving this open issue. 
We build upon the results presented in \cite{jan} and introduce an
optimization framework that is versatile and fast, eliminates
parameter dependencies, and ultimately produces highly accurate
rational filters with respect to a metric tightly bound to the quality
of the ideal filter. Numerical tests show that an eigensolver equipped
with our spectral projector converges with a rate that is practically
independent from the search space size, the number of poles of the rational
filter and the number of iterations required.

When $A$ is a Hermitian matrix, the corresponding rational filter is
real-valued and symmetric with respect to the mapping
$(x - x_0) \leftrightarrow (x_0 -x)$, where $x_0$ is the center point
of the interval $[a,b]$. Taking into consideration the complex
conjugation and parity symmetries, \wina\ write $r$
as a rational function of order  $(4m-1,4m)$,
\begin{equation}
\label{def:filters}
r(x) := \ratcp{r}{\beta,z}(x) := \sum_{i=1}^{m} \frac{\beta_i}{x - z_i} +
\frac{\conj{\beta_i}}{x - \conj{z_i}} - \frac{\beta_i}{x + z_i} -
\frac{\conj{\beta_i}}{x + \conj{z_i}}, \quad x \in \R,
\end{equation}
where $m \in \N$, $\beta=(\beta_1,\dots,\beta_m) \in \C^m$, and
$z=(z_1,\dots,z_m) \in (\HR)^m$, with $\HR$ being the right quadrant of the upper
half of $(\C \setminus \R)$ with origin in $x_0$.
With this setup, the
problem to be addressed is how to select, for a fixed degree
$m$\footnote{Strictly speaking the degree of $r$ is $4m$. In the rest
  of the paper we will stick to a more intuitive notion of
  degree which refers to the number of poles in $\HR$
  corresponding to the range of the index $i$ in
  \eqref{def:filters}}, the
coefficients $\beta_i$ and the poles $z_i$ such that the corresponding
rational function $r(x)$ approximates the ideal filter $\indi{a,b}$
according to a predetermined metric. Our aim is to build an
optimization framework and select an appropriate metric such that the
outcome is a filter $r$ stabilizing the convergence of the eigensolver.

Due to the discontinuity of the indicator function $\indi{a,b}$, 
the problem of determining the best coefficients and
poles for $r(x)$ is tackled
using a non-linear weighted least-squares approach. For a given
interval $[a, b]$, one aims to minimize the objective function
\begin{equation}
 \label{eqn:SLISE}
 f_{\wgtf}(\beta, z) := \int_{-\infty}^{\infty} \! \wgtf (x)
 \, (\indi{a,b}(x) - \ratcp{r}{\beta,z}(x))^2 \ \mathrm{d}x, \quad \text{where} \
 \beta \in \C^m, z \in (\HR)^m,
\end{equation}
over $\beta$ and $z$ for some fixed $m\in\N$ and a weight function
$\wgtf (x): \R \rightarrow [0,\infty)$, which is even and piecewise
constant. This optimization framework, termed \SLISE{} in \cite{jan},
provides a comprehensive parameterization of rational filters.  The
resulting \SLISE{} filters have proven to be competitive with previous
rational filters, such as \textsc{Gauss-Legendre}~\cite{feast} and
\textsc{Zolotarev}~\cite{guettel}.

The \SLISE{} framework is independent of the specific eigensolver in
which the function $r$ is plugged in and used as a spectral filter.
At glance, a filter optimized through this framework should
perform well independently from the target eigenproblem. In practice,
the effectiveness of a filter depends indirectly from the eigenvalue
distribution around the interval $[a, b]$ through the choice of the
weight function $\wgtf$.
In other words, despite its versatility, the \SLISE{
framework outputs filters whose quality is sensitive to the ad-hoc choice
of weight functions and the piece-wise intervals defining them: small
changes in the choice of $\wgtf (x)$ greatly influence the
effectiveness of the resulting filter.

\paragraph{Contributions}
Building on top of the \SLISE{} framework, this work addresses
problematic aspects of such optimization and ultimately provides a
solution to the open issue of how a spectral filters influence load
balancing. In detail, we identify a number of main contributions.  We
improve the performance of the unconstrained minimization process by
substituting the \textsc{Levenberg-Marquardt} with the
\textsc{Broyden-Fletcher-Goldfarb-Shanno} (\BFGS{}) algorithm
\cite[Chapter 6]{numericaloptmization}. Likewise, when \SLISE{} is
used in combination with box-constraints, it comes natural to extend
\BFGS{} to the \LBFGSB{} algorithm \cite{lbfgs1,lbfgs2,lbfgs3}.  Using
the \BFGS{} family of algorithms results in a substantial reduction of
time-to-solution, which in turn is a necessary requirement to reduce
the objective function residual and, at the same time, increases the
accuracy of the \SLISE{} filters. We increase the accuracy by
casting the problem of selecting $\beta_i$s and $z_i$s in terms of
finding the corresponding rational function $r(x)$
that 
minimizes the Worse-case Convergence Rate (\WCR{})\footnote{This
  metric is defined in the next section.}. The relevance of this
metric resides in the fact that the ideal filter $\indi{a,b}$ has the
lowest possible value for \WCR{}, which is $0$. 

In order to use the \WCR{} metric effectively, we embed the \SLISE{}
framework, equipped with the \BFGS{} algorithm, within a second
minimization process. This process has the explicit goal of minimizing
the \WCR{} metric with respect to the weight function $\wgtf$. We
attain this target by using the derivative-free \textsc{Nelder-Mead}
algorithm. The by-product of this process is eliminating the
dependence on the arbitrary choice of $\wgtf$ in the definition of the
objective function $f_\wgtf(\beta,z)$. The net result is a 
parameter-free minimization framework with an enhanced usability and
productivity. When used in interior eigensolvers based on subspace
iteration, we observe that the rational filters obtained with the new
minimization framework outperform state-of-the-art filters. The
convergence rate of the eigensolver becomes almost independent from
the size of the search subspace and the number of poles
used. Consequently, the eigensolver is more robust in terms of
convergence rate and does not require tweaking of the parameters
associated with the spectral projection. In turn, this enhanced
behavior of the eigensolver facilitates the load balancing when
executed on parallel platforms. We termed this enhanced minimization
framework, and the corresponding rational filters it produces,
\ESlise{}.

\paragraph{Related work}
The interpretation of spectral projectors as rational (filter)
functions of matrices was discussed
in~\cite{feast_subspace} for FEAST, and
in~\cite{Ikegami20101927,ikegami2010contour} for Sakurai-Sugiura-type
eigensolvers. Rational filters were also proposed early on for signal
processing by
Murakami~\cite{weko_67773_1,weko_69671_1,weko_28810_1,weko_28690_1,weko_18206_1}.
In recent years, filters have been treated as a parameter that can be
designed via optimization methods.  Van Barel~\cite{VanBarel2016346}
suggested a non-linear Least-Squares approach for non-Hermitian
filters to be used in the Sakurai-Sugiura framework, while Xi and
Saad~\cite{saad16LeastSquares} described linear Least-Squares
optimized filters for the Hermitian FEAST eigensolver. Van Barel's approach
is based on the discrete $\ell_2$ norm, not a functional approximation
approach, and does not support constraints optimization. Xi and Saad
present a linear Least-Squares minimization method where only the
coefficients of the rational function are optimized.  For FEAST,
G\"uttel et al.~presented a first approach to minimizing the \WCR{} in
\eqref{eqn:minwcr}. They derived a set of \textit{generalized}
\textsc{Gauss-Legendre} filters, parameterized by one variable only,
with respect to which they minimized the \WCR{} functional.  The
resulting \WCR{} values were smaller than for unparameterized
\textsc{Gauss-Legendre} filters, but not as small as for
\textsc{Zolotarev} filters
\cite{feast,guettel,guettel_correction,selfadjoint} which offered the best
\WCR{} so far.  This observation motivated a more rigorous
parameterization of a subset of rational filters, that is, \SLISE{}
filters~\cite{jan}, so as to benefit from a reduced number of parameters within
\WCR{} minimization.



\paragraph{Organization} The remainder of this paper is organized as follows.
In Sec.~\ref{sec:method} we introduce the reader to spectral filters and the
general mathematical setup. In Sec.~\ref{sec:slise}, we review the
\SLISE{} framework and introduce efficient imposition of
box-constraints on rational filters through \LBFGSB{}.  In
Sec.~\ref{sec:min_wcr}, we illustrate the minimization scheme to
reduce the \WCR{} of \SLISE{} filters, which in turns eliminates the
dependence on weight functions. In Sec.~\ref{sec:experiments}, we
present a set of numerical experiments comparing our new filters to
the state-of-the-art and illustrate their numerical properties and
advantages. The last section summarizes our results and provides a
perspective on their impact on the load balancing of parallel interior
eigenvalue solvers.


\section{Methodology}
\label{sec:method}

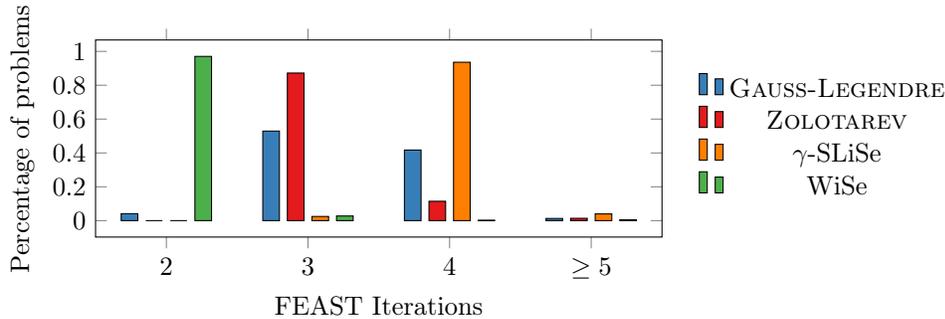
\begin{figure}
	\centering
	\tikzsetnextfilename{iterations_gamma_1.1_all}
	\pgfplotstableread[col sep = 
	comma]{figures/its_1.1.dat}\its
	\begin{tikzpicture}
	\begin{axis}[
	xmin=1.5,
	xmax=5.5,	
	xlabel=\FEAST{} Iterations,
	ylabel=Percentage of problems,
	ybar=3pt,
	xticklabels={2,...,4,$\ge 5$},xtick={2,...,5},
	ytick={0,0.2,...,1},
	bar width=0.017\textwidth,
	legend style={at={(1.05,0.5)},anchor=west},
	legend style={draw=none},
	width=0.7\linewidth,
	height=0.2\textheight,
	]	
	\addplot [fill=GAUSSC] table[x index={0},y expr=(\thisrowno{1}/2117)] 
	{\its};	
	
	\addplot [fill=ZOLOC] table[x index={0}, y expr=(\thisrowno{2}/2117)] 
	{\its};
	
	\addplot [fill=SLISEC3] table[x index={0}, y 
	expr=(\thisrowno{4}/2117)] 
	{\its};
	
	\addplot [fill=SLISEC] table[x index={0}, y expr=(\thisrowno{3}/2117)] 
	{\its};
	\legend{\textsc{Gauss-Legendre}, \textsc{Zolotarev}, 
	$\gamma$-\SLISE{}, \text{\ESlise{}}}
	\end{axis}
	\end{tikzpicture}
	\caption[Counts of iterations]{\FEAST{} iterations for 
	different filters with gap parameter $G=0.95$
	to solve $2117$ benchmark eigenproblems, for an eigencount multplier of 
	$C=1.1$.
	Our new \ESlise{} filter outperforms the others, i.e. 
	generalized \textsc{Gauss-Legendre} \cite{feast,guettel}, 
	\textsc{Zolotarev} \cite{guettel}, and $\gamma$-\SLISE{} \cite{jan}.
	Details are discussed in 
	Section~\ref{sec:experiments_benchmarkset}.
	}
	\label{fig:its_1.1}
\end{figure}


Contour-based eigensolvers were originally conceived for the
solution of the generalized interior eigenvalue problem
\begin{equation}
\label{eqn:eigenprob}
A v = \lambda B v, \quad \lambda \in [a,b],
\end{equation}
where $A,B\in \C^{n \times n}$ are Hermitian and $B$ is
positive definite, $v \in \C^{n}\setminus\set{0}$, $a<b$ and
$n \in \N$.
A spectral projector can be defined as the integral of the matrix
resolvent $(A - zB)^{-1}$ along a contour $\Gamma$ in the complex
plane $\C$ enclosing the interval $[a,b] \subset \R$. 
Without loss of generality, one can linearly map $[a,b]$ to the
standard interval $[-1,1]$ and select an integration contour around
it.
It is standard practice to compute the contour
integral via numerical quadrature (e.g. \textsc{Gauss-Legendre})
\begin{equation}
\label{eqn:countour}
r(A,B) := \sum_{i} \beta_i (A - Bz_i)^{-1} B\approx \frac{1}{2\pi i}
\oint_\Gamma \frac{\text{d}z}{A - Bz}B,
\end{equation}
with $\beta_i,\ z_i \in \C$.
When used in combination with a subspace iteration scheme, $r(A,B)$
projects a given set of vectors $Y$ onto an invariant subspace of the
spectrum corresponding to the eigenvalues within the interval $[-1,1]$
\cite{feast_conv}.  In practice, spectral projection exchanges the
direct solution of the eigenproblem for that of many independent
linear systems with multiple right-hand-sides
\begin{equation}
\label{eq:linsys}
(A - Bz_i) V = \beta_iB Y .
\end{equation}
Because each linear system can be solved independently from the
others, this class of eigensolvers naturally lends itself to
multiple layers of parallelism, making contour-based eigensolvers
especially well suited for today's increasingly parallel computer
architectures. As shown in several recent publications, the
performance of the eigensolver depends on the effectiveness of the
spectral filter $r(A,B)$
\cite{feast_conv,feast_subspace,guettel,guettel_correction,eigencount,ifeast}.
Recently, the authors of \cite{jan} proposed a numerical optimization approach alternative to the standard quadrature rules. By minimizing the objective function of
\eqref{eqn:SLISE}, they propose a new class of rational filters,
termed \SLISE{}, which perform better than the filters currently in
use, on a large number of representative eigenproblems. Despite such an
advance, the \SLISE{} framework showed a few shortcomings, such as
slow convergence and lack of efficient support for
box-constraints.
These box-constraints---defined as upper and lower bounds on the
imaginary parts of each $z_i$---can substantially influence the
time-to-solution in iterative linear system solvers. Having a 
time-to-solution comparable across all linear systems is a crucial element to
load-balance a parallel eigensolver based on spectral projection.

\paragraph{The \SLISE{} filters}
The \SLISE{} minimization framework aims to approximate the
indicator function $\ind$ by rational filters $r(x)$ of a fixed degree
$m$. This approximation is obtained by minimizing the
objective function $f_\wgtf(\beta,z)$ from \eqref{eqn:SLISE}. In the \SLISE{} framework, a new filter is obtained as follows: Given
a fixed weight function $\wgtf$ and an $m\in\N$, \SLISE{} takes an
existing rational filter $\ratcp{r}{\bar{\beta},\bar{z}}$, where
$\bar{\beta} \in \C^m $ and $\bar{z} \in (\HR)^m$, and derives a
new rational filter $\ratcp{r}{\hat{\beta},\hat{z}}$, such that
$(\hat{\beta},\hat{z})$ solves the minimization problem
\begin{equation}
 \label{eqn:minslise}
 \argmin{\beta \in \C^m, z \in (\HR)^m}{f_{\wgtf}(\beta,
  z)}.
\end{equation}
This minimization problem is non-linear and non-convex and therefore
difficult to solve due to the non-existence of closed-form solutions.
Yet, the objective function $f_{\wgtf}$, as well as its gradient
\begin{equation}
	\nabla f_{\wgtf} =
	(\nabla_{\beta_1} f_{\wgtf}, \dots, \nabla_{\beta_m} f_{\wgtf},
	\nabla_{z_1} f_{\wgtf}, \dots, \nabla_{z_m} f_{\wgtf})^\top,
\end{equation}
are differentiable and can be computed through a small number of matrix
operations \cite{jan}. In this setup, one can make use of a wide range of existing
numerical minimization methods.  \wina  obtain \SLISE{}
filters by employing two such minimization methods, gradient descent
and \textsc{Levenberg-Marquardt} (LM).

While LM makes for an effective minimization scheme, it may require
up to thousands of iterations to converge to a satisfactory value for
the residual level of $f_{\wgtf}$. Executing an efficient minimization
becomes a pressing problem in the case of box-constrained
optimization, when the LM algorithm cannot be used and gradient
descent requires up to millions of iterations to converge, which
translates in a significantly larger amount of computing time over the
unconstrained case.
In addition, and most importantly, the quality of
a resulting filter depends on the choice of weight function
$\wgtf$, which is not automatic and requires an experienced user to
follow a set of guidelines. In the following, we illustrate
a minimization scheme that ensures speed of convergence,
supports box-constraints, and eliminates the dependence on the custom
choice of weight functions $\wgtf$.

\paragraph{The new minimization scheme}
In the rest of the paper, we refer to $r$ as a rational
  filter, and, without loss of generality, consider only the case
$r(A,B=I)=r(A)$. As seen in the previous section, if $A$ is a
Hermitian matrix, the corresponding rational function $r(x)$ is forced
to be real and symmetric and can be expressed with a subset of
poles and coefficients as in \eqref{def:filters}.  Since the
minimization of the objective function in \eqref{eqn:SLISE} is
completely general, the resulting filter is independent of the
specific subspace iteration eigensolver and can be plugged in any eigensolver of
this type.
Nonetheless, for practical purposes, we use the
\FEAST{} eigensolver \cite{feastsolver} as a reference
algorithm. Given and exact value
$\lambda_j\in[-1,1]$, \FEAST{} computes an approximate
eigenpair $(\boldsymbol{q}_j,\hat{\lambda}_j)$ with a residual vector
norm equal to
$\norm{A \boldsymbol{q}_j - \hat{\lambda}_j \boldsymbol{q}_j}$. Such
residual converges linearly with a convergence rate given by
$\abbs{\gamma_{\rm out}/\gamma_{\rm in}}$, where $\gamma_{\rm out}$ 
($\gamma_{\rm in}$) is related to the maximum (minimum) value of
the filter outside (inside) a neighborhood enclosing the $[-1,1]$
interval%
~\cite[Theorem 5.2]{feast_subspace}. Consequently, the convergence rate
depends both on the spectrum of the given matrix $A$ and the
spectral filter of choice. 

Although the actual
convergence rate will vary for different spectra, a filter-dependent
upper bound is given by the Worst-case Convergence Rate
(\WCR{}). The \WCR{} applies to a variety of other
eigensolvers based on spectral projection such as the block Sakurai-Sugiura-Rayleigh-Ritz
method \cite{sakurai_cirr} and its non-iterative variant 
\cite{block_ss_rr_accuracy}.  As defined in \cite{guettel}, the
worst-case convergence rate satisfies the following theorem
\begin{theorem}[{\cite[Th.2.2]{guettel}}]
 \label{th:wcr}
 Given a rational filter $r$ and a fixed \textit{gap parameter} $G \in 
 (0,1)$, the \FEAST{} method converges linearly, with probability one,
 at a convergence rate no larger than
 \begin{equation}
  \label{eqn:worst}
  w_G(r) = \frac{\max_{x \in [-\infty,-G^{-1}] \cup [G^{-1},\infty]}\abs{r(x)}}{\min_{x \in [-G,G]}\abs{r(x)}},
 \end{equation}
 as long as no
 eigenvalues lie within $[-G^{-1},-G]\cup[G, G^{-1}]$.
 The occurring probability stems from choosing the initial subspace within
 the \FEAST{} method at random.
 \end{theorem}
 Since Theorem~\ref{th:wcr} implies that
 $\abbs{\gamma_{\rm out}/\gamma_{\rm in}} \le w_G(r)$, for an appropriate
 $G$, a smaller \WCR{} value $w_G(r)$ implies faster worst-case convergence. As we
 already mentioned in the introduction, minimizing \WCR{} for the
 \SLISE{} filters points out to which filters
 best approximate the ideal filter $\ind$ (that has indeed the optimal
 bound $w_G(\ind) = 0$). Based on the considerations above, we can now
 define the following optimization problem.
 \begin{definition}
 Given $G \in(0,1)$, $m\in\N$, and $\ratcp{r}{\beta,z}$ a rational
 filter as defined in \eqref{def:filters}, an \textit{optimal rational filter} is one
 solving the minimization problem
\begin{equation}
 \label{eqn:minwcr}
 \argmin{\beta \in \C^m, z \in (\HR)^m}{w_G(\ratcp{r}{\beta,z})}.
\end{equation}
\end{definition}


In general, the \WCR{} is a non-linear, derivative-free function.  Its
formulation makes it difficult to determine further mathematical
properties, such as convexity or continuity.  Conventional methods,
like steepest descent, cannot be applied. Additionally, in
derivative-free minimization, the number of $w_G$ function evaluations
may become intractable very quickly, even for a modest increase of the
rational filter degree $m$. These observations cause this
minimization problem to be especially challenging. Instead of solving
the problem as formulated in \eqref{eqn:minwcr}, we propose a modified
minimization problem that combines the existing \SLISE{} framework,
solving for $(\beta,z)$ while $\wgtf$ is fixed, as in
\eqref{eqn:minslise}, with minimization of the \WCR{} with respect to
$\wgtf$, seeking a better $\wgtf$ while $(\beta,z)$ is fixed.
These two minimization problems
\begin{eqnarray}
    \label{eqn:minWCR}
\left\{ 
\begin{array}[l]{l l l}
  \beta, z & \leftarrow  \argmin{\beta, z}{f_{\wgtf}(\beta,
  z)} & \textrm{for a fixed $\wgtf$} \\
 \wgtf  & \leftarrow 
 \argmin{\wgtf}{w_G(\ratslise{r}{\beta}{z}{\wgtf})}
 &  \textrm{for a fixed pair $(\beta,z)$}
\end{array}
\right.
  \end{eqnarray}
are clearly not independent.
The \WCR{} it is minimized solely with respect to the weight function
$\wgtf$---where we have indicated explicitly the dependence of the
rational filter on the weight function in \eqref{eqn:minWCR}---but it
is a non-linear, derivative-free function. As such, the WCR{} depends
on the whole $r$ which, in turn, depends on the minimization of the
objective function $f_\wgtf$. In other words, we now have to solve
two non-linearly dependent minimization problems, which need to be
solved self-consistently. We will see in Section \ref{sec:min_wcr} how
we implement this process in a nested loop fashion, where the \SLISE{}
process is executed within the \WCR{} minimization and convergence is
reached self-consistently by continuously swapping between the two
minimizations. Solving \eqref{eqn:minWCR} is now a tractable problem,
even if it calls for sophisticated algorithms, and is computationally
very intensive, requiring many repeated invocations of the \SLISE{}
minimization.

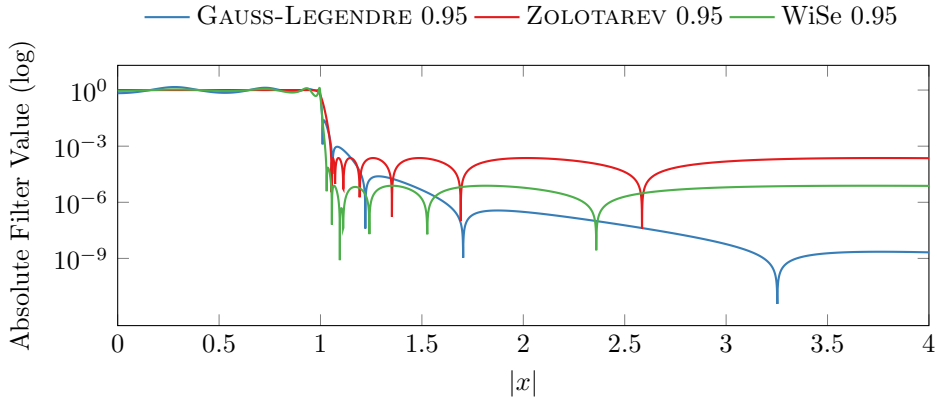
\begin{figure}
	\centering			
	\tikzsetnextfilename{filters_log}
	\begin{tikzpicture}
	\begin{axis}[
	xmin=0,
	xmax=4,
	ymode=log,
	ytick={0.000000001,0.000001,0.001,1},
	yticklabels={$10^{-9}$,$10^{-6}$,$10^{-3}$,$10^{0}$},
	xlabel={$\abs{x}$}, ylabel={Absolute Filter Value (log)},
	width=0.95\linewidth,
	height=0.24\textheight,
	legend style={at={(0.5,1.1)}, anchor=south,legend columns=-1},
	legend style={draw=none},
	]
	\addplot[color=GAUSSC, thick] table[x index={0},y 
	expr=abs(\thisrowno{1})] {figures/gl_4_95.dat};	 
	\addplot[color=ZOLOC, thick] table[x index={0},y 
	expr=abs(\thisrowno{1})] {figures/zolo_4_0.95.dat};	
	\addplot[color=SLISEC, thick] table[x index={0},y 
	expr=abs(\thisrowno{1})] {figures/wcr_4_95.dat};
	
	\addlegendentry{\textsc{Gauss-Legendre} $0.95$}
	\addlegendentry{\textsc{Zolotarev} $0.95$}
	\addlegendentry{\ESlise{} $0.95$}
	\end{axis}
	\end{tikzpicture}
	\caption{Logarithmic plot of different $16$-pole rational filters, showing 
	the state-of-the-art \textsc{Gauss-Legendre} and 
	\textsc{Zolotarev} rational filters \cite{guettel,feast}, 
	alongside our new 
	filter candidate, \ESlise{}. A moderate gap parameter 
	$G=0.95$ was chosen for \textsc{Zolotarev}, \ESlise{}, and 
	\textsc{Gauss-Legendre}.	
	Our new \ESlise{} filter provides the best \WCR{}, which can be seen 
	in offering 
	the sharpest slump for $x$ around $1$, and maintaining a constantly 
	close approximation of the indicator function $\ind$ across the whole domain of real numbers. 
    }
	\label{fig:filters_log}
\end{figure}

In order to increase the performance of the self-consistent
minimization, we introduce the \BFGS{} algorithm \cite[Chapter
6]{numericaloptmization} within the unconstrained \SLISE{}
minimization process. Similarly, for box-constraints, we present an
embedding of the \LBFGSB{} minimization method
\cite{lbfgs1,lbfgs2,lbfgs3} into \SLISE{}.  By formulating the \WCR{}
minimization as a nested process, we additionally solve the issue of
weight function selection, which is one of the open issues of
\SLISE{}.  The net result is an extension of the \SLISE{} framework
toward rational filters for faster convergence, without the need to
select the weight functions by hand. Our new rational filters, termed
\ESlise{}, outperform state-of-the-art filters (see
Figure~\ref{fig:its_1.1} for experimental results and
Figure~\ref{fig:filters_log} for filter plots). In particular, we prove
that \textsc{Zolotarev} filters do not provide best worst-case
convergence, despite their optimality in approximating the indicator
function with respect to the $\infty$-norm.


\section{Efficient computation of \ESlise{} filters}
\label{sec:slise}
In this section, we introduce the use the \BFGS{} algorithm, which
yields faster convergence and better box-constrained rational filters
than previous implementations. Moreover, the extended \LBFGSB{}
successfully addresses open issues that appear in box-constrained
filters~\cite{guettel,jan,ifeast}.

\subsection{Accelerating \SLISE{}}
\label{sec:bfgs}

When using the 
\BFGS{} algorithm 
to solve the minimization problem in \eqref{eqn:minslise}, we end up
reducing substantially the number of function evaluations needed to
reach convergence. Seemingly minor, this improvement is actually
essential for an effective embedding of \SLISE{} into a scheme that is
based on the minimization of the \WCR{}. The \BFGS{} algorithm belongs
to the class of quasi-\textsc{Newton} methods. It approximates a local
minimizer iteratively, in a manner similar to the popular
\textsc{Gauss-Newton} algorithm, which is
Hessian-based. However, unlike \textsc{Gauss-Newton}, \BFGS{} does not
require the exact Hessian $\nabla^2 f$, and uses an
approximation instead. The minimum requirement for the algorithm to
work is that the function $f$ has a quadratic expansion in
Taylor series near the minimum. Thanks to this weaker
condition, \BFGS{} guarantees convergence also for non-smooth and
non-convex functions.

The standard implementation of the \BFGS{}-variant in Algorithm~\ref{algo:bfgs}
does not offer support for real-valued objective functions of complex 
arguments, such as our $f_{\wgtf}$ from \eqref{eqn:SLISE}.
This problem can be overcome by a conversion of $f_{\wgtf}$ and $\nabla 
f_{\wgtf}$ to real arguments. In the
case of a generic function of a complex variable $g:\C^n \rightarrow
\R$, one can
separate the real from the imaginary parts \cite{complexoptimization}
and instead minimize the function $\tilde{g}:\R^{2n} \rightarrow\R$, defined as
\begin{subequations}
 \label{eqn:transformation}
 \begin{equation}
  \label{eqn:complex_transform}
  \tilde{g}(\colvecalt{a}{b}) := g(a + \iu b), \quad \text{for} \ a,b\in
  \R^n,
 \end{equation}
by computing descent directions from its gradient
 \begin{equation}
  \label{eqn:complex_transform2}
  \nabla \tilde{g} (\colvecalt{a}{b}) =
  \colvecalt{\operatorname{Re}\nabla g (a + \iu
  b)}{\operatorname{Im}\nabla g (a + \iu b)},
  \quad \text{for} \ a,b\in
  \R^n.
 \end{equation}
\end{subequations}
The same mapping can be applied to the \SLISE{} functional
$f_{\wgtf}: \C^m \times (\HR)^m \rightarrow \R$ because it operates on a
subset of $\C^{2m}$, where $m\in\N$ is the degree of the rational
filter. In this case, one can think of the complex vectors $\beta$
and $z$ as being part of a vector $v = (\beta\  z) ^\top$ and define
$\tilde{f}:\R^{4m} \rightarrow\R$ such that
 \begin{equation}
  \label{eqn:complex_transform3}
  \tilde{f}\colvecalt{\colvecalt{\Re(\beta^\top)}{\Re(z^\top)}}{\colvecalt{\Im(\beta^\top)}{\Im(z^\top)}}
  := f(\Re(\beta) + \iu \Im(\beta), \Re(z) + \iu \Im(z) ).
 \end{equation}

Starting at an initial point $x_0 = \left(\Re(\beta)\ \Re(z)\ \Im(\beta)\
  \Im(z)\right)^\top$, 
\BFGS{} computes iterates $x_k$ that converge to a local minimizer of
$\tilde{f}$ as $k \in \N$ increases, employing the descent directions
\begin{subequations}
 \begin{equation}
  \label{eqn:descent}
  p_k := - H_k \ \nabla \tilde{f}(x_k),
\end{equation}
and a line search which guarantees that the secant
condition is satisfied (see {\tt line \ref{algo:bfgs_linesearch}} of Algorithm
  \ref{algo:bfgs}).  $H_k$ is an approximation to the inverse
 Hessian of $\tilde{f}$ and is recursively defined as
 \begin{equation}
  H_0 := I_{4m}, \quad H_{k+1} := (I_{4m} - \frac{s_k y_k^{T}}{y_k^T
  s_k}) \ H_k \ (I_{4m} - \frac{y_k
  s_k^{T}}{y_k^T s_k}) +  \frac{s_k s_k^{T}}{y_k^T s_k},
 \end{equation}
 with
 \begin{equation}
  s_k := x_{k+1} - x_{k}, \quad y_k := \nabla \tilde{f}(x_{k+1}) -
  \nabla
  \tilde{f}(x_{k}),
 \end{equation}
 where $x_k, s_k, y_k, p_k \in \R^{4m}$ and $H_k \in \R^{4m \times 4m}$ for 
 some $m\in\N$.
\end{subequations}
The formulation through the \BFGS{} algorithm converges
faster than the previous minimization algorithms used by the \SLISE{}
framework (see Figure~\ref{fig:residuals} for the box-constrained case
that is discussed in the following subsection).
The conversion of the minimization functional to real-arguments 
allows one to use not only the \BFGS{} scheme, but also various other minimization
algorithms (such as those in the minimization algorithm collection
\tool{NLOpt} \cite{nlopt}). Despite such an advantage, most
alternatives do not yield any substantial improvements over \BFGS{}.

\begin{algorithm}[t]
 \caption{(Unconstrained \SLISE{} through \BFGS{} Algorithm).}
 \label{algo:bfgs}
 \begin{algorithmic}[1]
  \Procedure{\SLISE{}}{$\beta$, $z$, $\wgtf$} 
  \State Define $x \gets \left(\Re(\beta), \Re(z), \Im(\beta),
  \Im(z)\right)^\top \in \R^{4m}$
  \State Compute the objective function $\tilde{f}(x)$ 
  and
 	 its gradient $\nabla \tilde{f}(x)$
  \State $H \gets I_{4m}$
  \While{$\nabla \tilde{f}(x) > \eps$} \Comment{Default value $\eps=10^{-8}$}  
  \State $p \gets -H \ \nabla \tilde{f}(x)$ \Comment{Obtain descent
    direction}
  \State Choose an $\alpha \in \R^+$ to minimize
  $\tilde{f} ( x + \alpha p)$ over $\alpha$ \label{algo:bfgs_linesearch}
  \Comment{Ensuring $s^\top y >0$}
  \State $s \gets \alpha \, p$
  \State $w \gets x + s$ \label{algo:bfgs_newx}
  \State $y \gets \nabla \tilde{f}(w) - \nabla \tilde{f}(x)$
  \State $H \gets (I_{4m} - \frac{s y^\top}{y^\top
  	s}) \ H \ (I_{4m} - \frac{y s^\top}{y^\top s}) +  \frac{s 
  	s^\top}{y^\top s}$ 
  	\Comment{Approximate inverse Hessian}
  \State $x \gets w$  
  \EndWhile
  \State $\beta' \gets (x_{1:m} + \iu x_{2m+1:3m})^\top, \quad z' \gets  (x_{m+1:2m} + \iu x_{3m+1:4m})^\top$ 
  \State \Return $\ratcp{r}{\beta',z'}$
  \EndProcedure
 \end{algorithmic}
\end{algorithm}


\subsection{Imposing box-constraints efficiently}

As described at the beginning of Section \ref{sec:method}, the
spectral projection at the base of the \tool{FEAST} eigensolver leads to the
solution of several independent linear systems with multiple RHS (see
\eqref{eq:linsys}).
In the case of very large and sparse systems, the use of direct
solvers is not feasible due to memory requirements.  In this case,
iterative solvers, such as \tool{GMRES} or \tool{CG}, are the natural
choice. For these
methods, time-to-solution and accuracy depend substantially on the
condition number of the resolvent matrices $(A - z_iI)$. When $A$ is
Hermitian, such condition number is, up to a constant factor, equal to
\begin{equation}
	\kappa(A- z_iI) = \frac{\max_{\lambda_a \in \sigma(A)}\abbs{\lambda_a - z_i}}{\min_{\lambda_b \in \sigma(A)}\abbs{\lambda_b
			- z_i}}.
\end{equation}
Since the filter is built to approximate the indicator function
$\ind$, the numerator of this equation is bound from above by
$(\max_{\lambda_a \in \sigma(A)}\abbs{\lambda_a} + 1)$, while the denominator is bound from below by
$\abbs{\Im(z_i)}$. Consequently, if the poles of the rational function
$\ratcp{r}{\beta,z}$ are close to the real axis, the condition number of some of the
resolvent matrices can be quite high. This
consideration motivated the introduction of the box-constraints
$\abbs{\Im(z_i)} \geq {\tt lb} > 0$ ($i=1, \ldots, m$) to the \SLISE{}
minimization process, where ${\tt lb}$ is a positive constant
representing the minimum distance of any pole from the real axis.

In the \BFGS{} algorithm, box-constraints can be included by projecting
the search direction onto the constraints. This is accomplished through the simple gradient projection
$\mathcal{P}(x - t \nabla\tilde{f}(x))$, where $t>0$, followed by a \BFGS{}
update treating the bounded components of $x$ as equality
constraints. In our case, the operator $\mathcal{P}: \R^{4m}
\rightarrow \R^{4m}$ projects only the imaginary
part of the poles $\abbs{\Im(z_i)}$ and takes consequently the following
form when acting on a vector $y\in\R^{4m}$
\begin{equation}
	\mathcal P(y)_j \coloneqq
	\begin{cases}
	\sign(y_j) \cdot \text{\texttt{lb}}, &
	\textrm{if}\ \abbs{y_j} >
	\text{\texttt{lb}} \ \text{and} \ j \in \set{ 3m+1,\ldots,4m },\\
	y_j & \text{otherwise, }
	\end{cases}
\end{equation}
for $j=1,\dots,4m$.

This approach is encoded in the \LBFGSB{} algorithm, which extends
projected gradient descent to the Hessian approximations
from \BFGS{}, and can be used to realize box-constrained
minimization efficiently in \SLISE{}, similarly to what is done in Algorithm
\ref{algo:bfgs}.
The \LBFGSB{} algorithm
has shown to converge quickly in our experiments, when compared with projected
gradient descent.
In terms of both speed and accuracy, the use of \LBFGSB{} places the
constrained \SLISE{} method on par with the unconstrained \BFGS{} algorithm.  
To illustrate the increase in performance caused by \LBFGSB{}, we
compare box-constrained minimization through our \LBFGSB{}
implementation against the projected gradient descent implemented in
the original \SLISE{}
framework. Figures~\ref{fig:projected_descent_residuals} and
\ref{fig:lbfgsb_residuals} show the number of function evaluations
carried out by the projected gradient descent and the \LBFGSB{}
algorithms, respectively. \LBFGSB{} requires four order of magnitudes
fewer evaluations than projected gradient descent, and converges to a
smaller residual.

So far, the procedure used to obtain the \SLISE{} filters depends on
the specific form of a given weight function $\wgtf$.  For some of
such weight functions, the outputted filters were shown to outperform
state-of-the-art rational filters. Yet, the only criterion known to
determine suitable weight functions is by comparing hand-crafted weight
functions on a large set of representative interior eigenproblems.
While guidelines for the construction of $\wgtf$ have been devised,
the choice of weight functions remains a complex issue. In
the following section, we propose an algorithm to obtain weight
functions, which yield \SLISE{} filters with reduced \WCR{}, and
overcome the necessity of selecting weight functions manually.

\begin{figure}
	\centering
	\subfloat[Projected gradient descent]{
		\label{fig:projected_descent_residuals}
		\tikzsetnextfilename{projected_descent_residuals}
		\begin{tikzpicture}
			\begin{axis}[
			scaled x ticks=false, 
			xtick={4000000, 8000000},
			xlabel={Function evaluations}, ylabel={
			$f_{\wgtf}$ (log)},
			enlargelimits=false,
			ymode=log,
			width=0.47\textwidth,
			height=0.2\textheight,
			]
			\addplot[color=GAUSSC, thick] table[x index={2},y index={1}] 
			{figures/box_residuals_projected_descent.dat};
			\end{axis}
		\end{tikzpicture}}
	\subfloat[\LBFGSB{}]{
		\label{fig:lbfgsb_residuals}
		\tikzsetnextfilename{lbfgsb_residuals}
		\begin{tikzpicture}
		\begin{axis}[
		xlabel={Function evaluations}, ylabel={$f_{\wgtf}$ (log)},
		enlargelimits=false,
		ymode=log,
		width=0.47\textwidth,
		height=0.2\textheight,
		]
		
		\addplot[color=ZOLOC, thick] table[x index={2},y index={1}] 
		{figures/box_residuals_lbfgsb.dat};
		\end{axis}
		\end{tikzpicture}}
              \caption[Residuals in box-constraints]{ Box-contrained
                minimization of the functional $f_{\wgtf}$ using
                projected gradient descent and \LBFGSB{} respectively.
                The setup is taken from the original publication
                \cite{jan}, using the $16$-degree \textsc{Zolotarev}
                filter as starting point and a lower bound of
                $\text{\texttt{lb}}=0.0022$ on the absolute value of
                the imaginary parts of the poles. The \LBFGSB{} method
                settles at a smaller residual and converges
                substantially faster, requiring only a few function
                evaluations only.  }
	\label{fig:residuals}
      \end{figure}
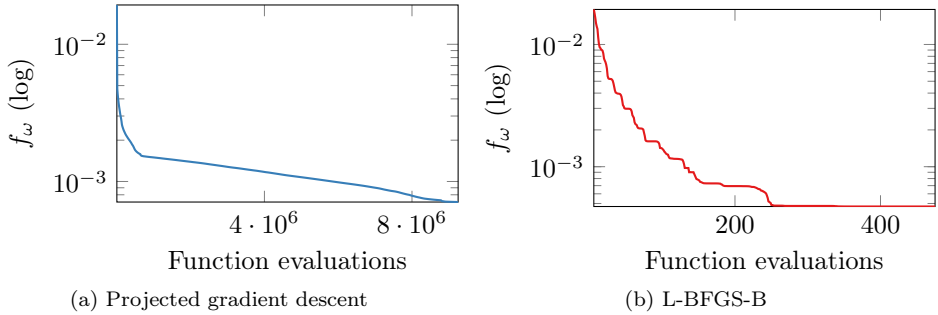
      


\section{\SLISE{} filters with reduced \WCR{}}
\label{sec:min_wcr}
In this section, we illustrate how to reduce the \WCR{} of a given 
\SLISE{} filter by improving on the choice of weight function $\wgtf$.
We achieve this by minimizing a new objective function, closely related to 
the \WCR{} of rational filters.

\subsection{Parameterization of weight functions}
\label{sec:param_weights}

Weight functions are even, non-negative, piecewise constant functions
that are used in the definition of the \SLISE{} functional $f_{\wgtf}$
in \eqref{eqn:SLISE}.  This means that a weight function can be
characterized in terms of $n\in\N$ intervals
$[x_i, x_{i+1}) \subseteq [0,\infty]$ and corresponding function
values $\wgtf(x\in [x_i, x_{i+1}))=\wgtf_i$, called weights,
where $\wgtf_i \in [0,\infty)$ for $i = 1,2,\dots,n$.  In their original
contribution, \wina{} obtained weight functions for the \SLISE{}
framework by following three guiding principles, derived from
experience: (i) gradual decrease in weights outside the search
interval $[-1,1]$, (ii) sufficient magnitude of weights inside
$[-1,1]$, and (iii) symmetry in weights about the interval endpoints
of $[-1,1]$.  While \SLISE{} filters following these guidelines could
outperform state-of-the-art \textsc{Gauss-Legendre} and
\textsc{Zolotarev} filters, some manual adjustment based on experience
remained necessary.

\begin{table}[t]
	\centering
	\begin{tabularx}{\linewidth}{ X  D  D D  D  D}
		\toprule
		$|x|\in$ & $[0, 0.95)$    & $[0.95, 1.05)$ & 
		$[1.05, 1.4)$ & $[1.4, 5)$ & $[5, \infty)$\\
		\midrule
		$\wgtfgamma(x)$ & 1 & 0.01 & 10 & 20 & 
		0                             \\
		\bottomrule
	\end{tabularx}
	\vspace{0.2cm}
	\caption{The $\wgtfgamma$ weight function
      }
	\label{tb:cuswei}	
\end{table}

An example of a weight function $\wgtfgamma$, which yields \SLISE{}
filters outperforming \textsc{Gauss-Legendre} filters, is given in
Table~\ref{tb:cuswei}.  This choice of weights suggests a natural way
of parameterizing weights and interval boundaries, so they can be
treated without distinction. For this purpose, we introduce a set
$\wgtfs$ of parameter vectors $v=(v_1, \dots, v_{2s-3})$, where $s \ge 2$
 equals the number of intervals to the right
of $0$ 
\begin{equation}
	V_s = \set{v \in [0,\infty)^{2s-3} \mid G < v_1 < 1 < G^{-1} < 
	v_2 < v_3< \ldots < v_{s-1}},
\end{equation}
for some gap parameter $G \in (0,1)$. 

A generic set of $v_i\in \wgtfs$ induces a weight function
$\wgtfsearch{j}$ with $j = 1,2,\dots,s$.
Following this parameterization,
Table~\ref{tb:cuswei} is rewritten as Table~\ref{tb:genwei}. The
parameters $v_1, v_2$ enclose $\pm 1$, but do not necessarily match
the endpoints of the gap $[G, G^{-1}]$.  The parameters $v_3,v_4$
reflect some more intervals of the weight function.  The remaining
parameters $v_5, \dots, v_7$ denote non-negative weights.  The weight
for the interval $[0,v_1)$ is fixed to $1$, as weight functions are
invariant under scaling within \SLISE{}. In this notation, the
$\wgtfgamma$ weight function from Table~\ref{tb:cuswei} translates
into the vector $(0.95, 1.05, 1.4, 5, 0.01, 10, 20) \in V_5$. As we
are going to illustrate in the next section, this parameterization
scheme allows for a systematic improvement of weight functions,
alongside a choice of weights and interval
endings.
\begin{table}
	\begin{tabularx}{\textwidth}{ X  D  D  D  D  D}
		\toprule
		$|x|\in$ & $[0, v_1)$ & $[v_1, v_2)$ & $[v_2, v_3)$
		& $[v_3, v_4)$ & $[v_4, \infty)$ \\
		\midrule
		$\wgtfsearch{j}(x)$ & $1$ & $v_5$ & 
		$v_6$ & $v_7$ & 0 \\
		\bottomrule
	\end{tabularx}
	\vspace{0.2cm}
	\caption{Parameterized weight function $\wgtf$, for 
	$s=5$.} 
	\label{tb:genwei} 
\end{table}

\subsection{Minimization of parameterized weight functions}

In order to compare the influence of distinct weight functions
$\wgtfsearch{j} \in \wgtfs$ on the minimization of
\WCR{}, we introduce a new objective
function
\begin{equation}
		\label{eqn:weightmetric}
		h := h_{\beta,z}(v) :=
		w_G(\ratslise{r}{\beta}{z}{v}), \quad \text{for} 
		\ v\in \wgtfs,
\end{equation}
where $\ratslise{r}{\beta}{z}{v}$ is computed by the
Algorithm~\ref{algo:bfgs}.
$h_{\beta,z}(v)$ is a functional of a given filter
$\ratcp{r} {\beta,z}$,
and it associates the vector of parameters $v$ with the \WCR{} of the corresponding
filter. As such, $h$ establishes a meaningful 
metric to quantify the performance of a weight function
$\wgtf$
. The minimization
\begin{equation}
\label{eqn:minweight}
\argmin{v \in \wgtfs}{h_{\beta,z}(v)},
\end{equation}
facilitates a systematic search for better weight functions and,
consequently, rational filters with smaller worst-case convergence.
For a given $\ratcp{r}{\beta,z}$, this is a non-linear,
derivative-free minimization problem, depending on only ($2s-3$)
variables.
However, for each optimization step
involving changes of $v$, the rational filter
$\ratslise{r}{\beta}{z}{v}$ has to be computed again by executing a
call to the Algorithm~\ref{algo:bfgs}.  The end result is a nested
optimization problem \eqref{eqn:minWCR} with, possibly, thousands of
calls to Algorithm~\ref{algo:bfgs}.

\begin{algorithm} [t] 
	\caption{(Local \WCR{} minimization)}
	\label{alg:min_weight_worst}
	\begin{algorithmic}[1]
		\Procedure{ReduceFilterWCR}{$v,\beta,z,G$}
		\State $G \gets \sqrt{G}$ \label{lln:shifting1} \Comment{For 
		shifting of filter}
		\State $n \gets \Call{Length}{v}$
		\State $w \gets 
		h_{\beta,z}(v)$ 
		\While{$\textrm{Res}(h) > 10^{-9}$} \label{alg:minww_res}
		\For{$i=1,2,\dots,n$} \label{lln:restart1} \Comment{Coordinate descent}
                \State{$\hat{v} \gets$}
                \Call{AdaptiveDifferentialEvolution}{$h_{\beta,z}(v(v_i))$}
                \If{$h_{\beta,z}(\hat{v}) < w$}
               		\State $v \gets \hat{v}, \ \textrm{and}\ w \gets h_{\beta,z}(\hat{v})$
                \EndIf
		\EndFor \label{lln:restart2}
                \State $v'\gets$ \Call{Nelder-Mead}{$h_{\beta,z}(v(v_3, 
		\dots, v_{n}))$} \label{lln:neldermead}
		\State $v' \gets (\abbs{v'_1}, \dots, 
			\abbs{v'_{n}}) \in \wgtfs$
			\label{lln:fixconstraints}
                \State{$\ratcp{r}{\beta',z'}\gets$
                  \Call{\SLISEtext{}}{$\beta,z,\wgtf'$}
                  \Comment{See Algorithm \ref{algo:bfgs}}}
		\State $w' \gets h_{\beta',z'}(v')$ 
		\Comment{New \WCR{} value}
		\State $\textrm{Res}(h) \gets \abbs{w - w'} \, / \, w$
                \Comment{Compute \WCR{} residual}
                \State $w \gets w', v \gets v', \quad \beta \gets \beta', \quad z \gets z'$
		\EndWhile
                \State $(\beta, z) \gets (G^{-2} \, \beta, G^{-2} \, z)$
                \label{lln:shifting2} \Comment{Shifting coefficients
                  and poles}
                \State \Return $\ratcp{r}{\beta,z}$
		\EndProcedure
	\end{algorithmic}
\end{algorithm}

Since the WCR functional cannot be expressed as a continuous function
of $v$, to solve \eqref{eqn:minweight}, we resort to using the
\textsc{Nelder-Mead} algorithm, a
prominent local,
derivative-free minimization schemes.\footnote{Local means that the algorithm starts at an
  existing point, at best, close to the sought after minimum.}
In our case,
\textsc{Nelder-Mead} generates competitive solutions quickly, but
suffers from stagnation at non-optimal points
\cite{neldermead_nonoptimal1,neldermead_nonoptimal2}.  To overcome
stagnation, we follow Carl Kelley's suggestion of restarting
\textsc{Nelder-Mead} at the current iterate with adjusted parameters
\cite{neldermead_nonoptimal}.  An explicit such parameter choice
exists only for the case of smooth functions, introduced as
\textit{oriented restart}. Since our functional $h$ is non-smooth, we
obtain a new parameter choice by perturbing the current iterate
carefully through \textit{coordinate descent}. Coordinate descent is a
simple, local, derivative-free minimization method, that performs
subsequent line searches along the coordinate directions, given some
starting point $v\in \wgtfs$. Independently of having detected
stagnation, we use coordinate descent systematically to obtain a new
starting point for each \textsc{Nelder-Mead} call.

The general scheme outlined above is described in
Algorithm~\ref{alg:min_weight_worst} and implemented using the
\tool{Julia} programming language.\footnote{The code is freely available at
  \url{https://github.com/SimLabQuantumMaterials/SLiSeFilters.jl}}
Given a weight function
$v \in \wgtfs$ and a filter $\ratcp{r}{\beta,z}$, the algorithm
chooses better weight functions from $\wgtfs$
iteratively. 
At each iteration of the {\bf while} loop, the weights of the current
filter are updated to reduce the \WCR{}, and from these, a new filter
is computed.
When the residual of the $h$ function $\textrm{Res}(h)$ (i.e. the relative difference of two subsequent $h$ values) falls below an established threshold
tolerance (see {\tt line \ref{alg:minww_res}} of
Algorithm~\ref{alg:min_weight_worst}), the algorithm returns the \SLISE{}
filter of the last iteration, which minimizes the \WCR{} among the
weight functions in the search space $\wgtfs$.
Each {\bf while} loop iteration follows three consecutive steps: (i) coordinate
descent, (ii) \textsc{Nelder-Mead}, (iii) computation of new \SLISE{}
filter and convergence check.
\paragraph{Coordinate descent} This is performed in
{\tt Lines~\ref{lln:restart1}~-~\ref{lln:restart2}}, improving the
coordinates of the parameter vector $v \in \wgtfs$ through a
separate minimization problem for each variable
$v_i$ and $i \leq 2s-3$, 
\begin{subequations}
	\begin{equation}
		\label{eqn:nelder_restart}
		\argmin{c\in I_i }{h_{\beta,z}(v(c))}, \quad 
		\text{where} \
		v(c) := (v_1, \dots v_{i-1}, \, c \, , v_{i+1} 
		\dots, 	v_{2s-3}) \in \wgtfs,
	\end{equation}
	while restricting the search space to a neighborhood $I_i$ of 
	$v_i$. For instance, for $s \ge 5$ the intervals $I_i$
        used are
	\begin{equation}
		I_i := \left\{\begin{alignedat}{3}
		&\ [ G, && 1],   && \quad \text{if } i=1, \\
		&\ [ 1, && G^{-1}], && \quad \text{if }   i=2 \\
		&\ [ G^{-1}, && v_{i+1}],  && \quad \text{if } i=3, \\
		&\ [ v_{i-1}, &&  v_{i+1}],  && \quad \text{if } 4 \le 
		i 
		<s-1, \\
		&\ [\, v_{i-1},  &&  3 \, v_i],  && \quad \text{if } 
		i=s-1, 
		\\
		&\ [ 0.1 \, v_i,   && 10 \, v_i],  && \quad \text{if } 
		s 
		\le 
		i \le 
		2s-3.
		\end{alignedat}\right.
	\end{equation}
\end{subequations}
We implement the coordinate descent minimization through the global,
derivative-free minimization scheme Adaptive Differential
  Evolution from the \tool{Julia} library \href{https://github.com/robertfeldt/BlackBoxOptim.jl}{BlackBoxOptim.jl}
\cite{blackboxoptim}.  Global minimization algorithms aim to
find the global minimizer within a region $I_i$, instead of converging
to a local minimizer starting from a given point. If the value of
\WCR{} has decreased, the solution of the coordinate descent
minimization $\hat{v}_i$ is used instead of $v_i$ for the successive
steps of the {\bf while} loop iteration. This step is executed to
prevent stagnation in the execution of the \textsc{Neldear-Mead}
minimization.

\paragraph{Nelder-Mead} In {\tt Line~\ref{lln:neldermead}},
\textsc{Nelder-Mead} is applied to a slightly modified version of the
minimization problem formulated in \eqref{eqn:minweight}: We keep
$v_1$, $v_2$ fixed at the values obtained from the coordinate descent
step. This choice is motivated by the high sensitivity of the
functional $h$ to changes in these variables, being
close to the endpoints of the gap $[G,G^{-1}]$. In practice,
\textsc{Nelder-Mead} is used to solve the 
minimization problem
\begin{subequations}
	\begin{equation}
	\label{eqn:neldermead_minimize}
	\argmin{ b_1, \dots, b_{r} \in \R }{h_{\beta,z}(v(b_1, 
		\dots, b_{r}))},
	\end{equation}
	where
	\begin{equation}
	v(b_1, \dots, b_{r}) := (v_1, v_2, b_1, \dots, 
	b_{r}) \quad {\rm and} \quad r=2s-5.
	\end{equation}
\end{subequations}
Because the classical \textsc{Nelder-Mead} algorithm performs 
unconstrained minimization only, we also
have to ensure that the minimizer $v'$ 
lies within the admissible intervals defined by $V_s$.
A straightforward approach would be to use a modified objective 
function within \textsc{Nelder-Mead}, defined as
\begin{equation}
\label{eqn:penalty}
\hat{h}(v) := \begin{cases} 
h(v), & \text{if } v\in 
\wgtfs, \\
\infty, & \text{otherwise,}\end{cases}
\end{equation}
to penalize invalid weight functions $v \notin \wgtfs$.  We verified
that this approach works, but our experiments have shown that it slows
down the convergence to the minimum. Based on our tests, we observed
that, in practice, most of the violations $v \notin \wgtfs$
stem from the selection of slightly negative weights by \textsc{Nelder-Mead}. This is likely caused by rounding errors, which may lead to a slight decrease of the \WCR{}
value for very small but negative $\wgtfsearch{j}$. Every other
violation of constraints seems to cause the opposite
of a reduction in \WCR{} value and is thus not chosen by
\textsc{Nelder-Mead}. To overcome this problem, we adopted a very
simple solution: 
We map the resulting minimizer of a \textsc{Nelder-Mead} into
$\wgtfs$ explicitly, by taking the absolute values
$(\abbs{v_1^{(k+1)}}, \dots, \abbs{v_{2s-3}^{(k+1)}})$ instead of a possible
negative (invalid) iterate $v^{(k+1)}$ (see \texttt{line}~\ref{lln:fixconstraints}).
While this approach only avoids violations of constraints caused by
the choice of  negative
weights, we did not
experience other violations of constraints 
in the minimizers returned by \textsc{Nelder-Mead}. For the sake of
completeness, we have experimented with other 
derivative-free minimization algorithms (those from the minimization
algorithm collection \tool{NLOpt} \cite{nlopt}, including a
constrained version of \textsc{Nelder-Mead}), none of which led to a
more competitive reduction of \WCR{} value in the same standard
setups.
 

\paragraph{\SLISE{} filter and convergence} As already mentioned, the \SLISE{} procedure is called multiple times within both
the \textsc{Nelder-Mead} and the
\textsc{AdaptiveDifferentialEvolution} procedures. The last call of
\SLISE{} is executed so as to calculate the residual of the \WCR{}
functional and check for convergence. If convergence is reached, the
algorithm returns values for the \WCR{}, the poles $z$ and the
coefficients $\beta$. The latter are rescaled, as it is explained in
the following Subsection~\ref{sec:shifting}, by introducing a
linear scaling transformation to improve the behavior of the resulting
filters at the interval $[-1,1]$ boundaries.

\subsection{Scaling the filter} 
\label{sec:shifting}

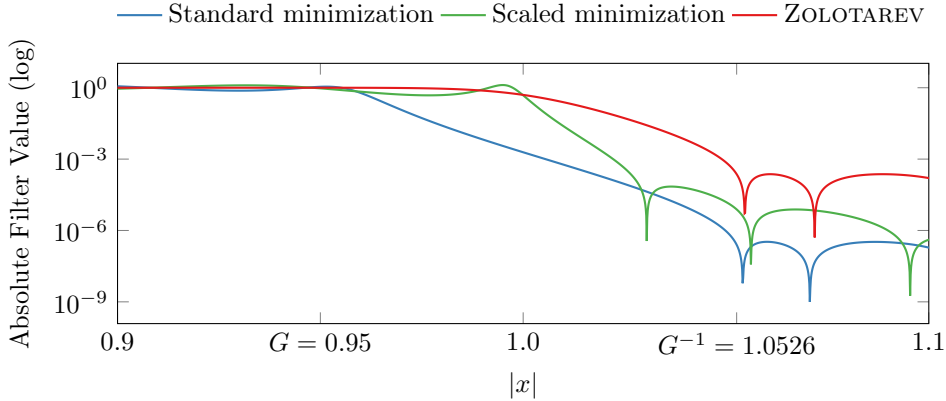
\begin{figure}
	\centering			
	\tikzsetnextfilename{filters_shifting_problem}
	\begin{tikzpicture}
	\begin{axis}[
	xmin=0.9,
	xmax=1.1,
	xtick={0.9,0.95,1,1.052631579,1.1},
	xticklabels={$0.9$,$G=0.95$,$1.0$,$G^{-1}=1.0526$,$1.1$},
	ymode=log,
	ytick={0.000000001,0.000001,0.001,1},
	yticklabels={$10^{-9}$,$10^{-6}$,$10^{-3}$,$10^{0}$},
	xlabel={$\abs{x}$}, ylabel={Absolute Filter Value (log)},
	width=0.95\linewidth,
	height=0.24\textheight,	
	legend style={at={(0.5,1.1)}, anchor=south,legend columns=-1},
	legend style={draw=none},
	]
	
	\addplot[color=GAUSSC, thick] table[x index={0},y 
	expr=abs(\thisrowno{3})] {figures/wcr_4_95_accuracy.dat}; 
	\addlegendentry{Standard minimization}
	
	\addplot[color=SLISEC, thick] table[x index={0},y 
	expr=abs(\thisrowno{1})] {figures/wcr_4_95_accuracy.dat};
	\addlegendentry{Scaled minimization}	
	
	\addplot[color=ZOLOC, thick] table[x index={0},y 
	expr=abs(\thisrowno{1})] {figures/zolo_4_0.95_accuracy.dat};
	\addlegendentry{\textsc{Zolotarev}} 
	
	\end{axis}
	\end{tikzpicture}
	\caption{Logarithmic plot of $16$-pole rational filters	for gap 
	parameter $G=0.95$.
	The standard filter is a solution of the original minimization 
	problem as of \eqref{eqn:minweight},
	while the scaled filter solves the new minimization problem in
	\eqref{eqn:newwcr}.
	Even though the scaling causes the WCR to
          increase from $\expnumber{4.95}{-6}$ to 
	$\expnumber{1.04}{-5}$ , this is still below the WCR of
        \textsc{Zolotarev} filter ($\expnumber{2.32}{-4}$).
	On the other hand, near $1$, the standard filter attains a function value of 
	about $10^{-3}$. 
	Hence, an eigenvalue inside the gap $[G,1]$
	could roughly half the actual convergence rate for the standard 
	filter.
	The other two filters are not affected by this problem.}
	\label{fig:filters_shifting_problem}
\end{figure}


Compared to previous filters such as \textsc{Gauss-Legendre} or
\textsc{Zolotarev}, our new \SLISE{} filters obtained from the
minimization problem \eqref{eqn:minweight} reduce the \WCR{} by up to
multiple orders of magnitude.  So far, we assumed to know the
appropriate gap parameter $G$ in advance and minimized our filters
accordingly.  This assumption may be too optimistic for some
eigenproblems: For a given $G$, we assume that no eigenvalues lie
within the interval $[-G^{-1},-G] \cup [G, G^{-1}]$.  Since, in
practice, this assumption may be violated, a rational filter with
small function value within $[G,1]$ may lead a slower convergence rate
than the \WCR{}, or no convergence at all.  An illustration of this
problem is given in Figure~\ref{fig:filters_shifting_problem}, and its
caption. Eigenvalues within $[1,G^{-1}]$ are less problematic because
they are not sought after by the eigensolver (for an in-depth
discussion, see Section~\ref{sec:feastexperiments}).  While
maintaining a competitive \WCR{} value, the issue described above can
be overcome by solving the slightly modified minimization problem
\begin{equation}
	\label{eqn:newwcr}
	\argmin{\beta \in \C^m, z \in (\C 
		\setminus \R)^m}{\frac{\max_{x \in [-\infty,-G^{-1}] \cup 
				[G^{-1},\infty]}\abs{\ratcp{r}{\beta,z}(x)}}{\min_{x \in 
				[-1,1]}\abs{\ratcp{r}{\beta,z}(x)}}}
\end{equation}
instead of minimizing $w_G$ as in \eqref{eqn:minwcr}.  For any
solution $(\beta,z)$ of this modified problem, the rational filter
$\ratcp{r}{\beta,z}$ offers a larger function value inside the entire
search interval $[-1,1]$ than outside in
$[-\infty,-G^{-1}] \cup[G^{-1},\infty]$.  When our filters are used in
iterative eigensolvers based on spectral projection, this behavior
ensures reliable and fast convergence.

Instead of modifying the \WCR{} minimization procedure,
we solve for the modified \WCR{} in \eqref{eqn:newwcr} by
introducing a linear scaling transformation $u(x) := \sqrt{G} \, x$.
We have
	\begin{align*}
	\label{eqn:minwcr_old}
	\argmin{\beta \in \C^m, z \in (\HR)^m}{w_{\sqrt{G}}(\ratcp{r}{\beta,z})}
	= &\argmin{\beta \in \C^m, z \in (\HR)^m}{\frac{\max_{x 
	\in [-\infty,-G^{-1}] \cup 
			[G^{-1},\infty]}\abs{\ratcp{r}{\beta,z}(u(x))}}{\min_{x \in 
			[-1,1]}\abs{\ratcp{r}{\beta,z}(u(x))}}}.
	\end{align*}
Since the function composition $\ratcp{r}{\beta,z} \circ u$ is a rational 
filter itself, if $\ratcp{r}{\beta,z}$ solves
\eqref{eqn:minwcr_old}, then
$\ratcp{r}{\beta,z} \circ u$ solves
\eqref{eqn:newwcr}.
Additionally, it follows from the 
definition of rational filters in \eqref{def:filters} that $\ratcp{r}{\beta,z} \circ u =      
\ratcp{r}{\sqrt{G^{-1}}\,\beta,\sqrt{G^{-1}}\,z}$ which characterizes 
the parameters of the resulting filter.
This scaling of rational filters through linear transformation 
is
incorporated in Algorithm~\ref{alg:min_weight_worst} in 
\texttt{lines~\ref{lln:shifting1}} and \texttt{\ref{lln:shifting2}}.

In our implementation of Algorithm~\ref{alg:min_weight_worst},
we used as initial parameters the weight functions
$\wgtfgamma$ used in \cite{jan} opportunely rescaled. For instance, for
$s=5$, we selected
\begin{equation*}
	v^{(0)}=({\sqrt{G}},{\sqrt{G^{-1}}},1.4,5,.01,10,20)
\end{equation*}
as the initial parameters characterizing the weight function,
for some gap parameter $G \in (0,1)$.
Accordingly, we chose both \textsc{Zolotarev} and \textsc{Gauss-Legendre} 
filters as initial conditions for \WCR{} minimization.
We obtained a number of new spectral filters, to which we refer as
\textbf{W}orst-Case Opt\textbf{i}mized Least-\textbf{S}quar\textbf{e}s
(\ESlise{}).  In the following section, we provide a number of
experimental tests illustrating the performance of these \ESlise{}
filters.


\section{Experiments}
\label{sec:experiments}

In this section, we compare the \textsc{Zolotarev} and generalized
\textsc{Gauss-Legendre} filters to the \ESlise{} filters from the
previous section in two different scenarios. First, in order to
inspect the worst-case performance, we compute the \WCR{} values for
different gap parameters $G$ and poles per quadrant $m$.  Second, we
use the filters in the \FEAST{} package and assess their performance
on two eigenproblems used in past literature \cite{feast,jan,guettel}.

We provide the \texttt{Julia}
library 
\href{https://github.com/SimLabQuantumMaterials/SLiSeFilters.jl}{SLiSeFilters.jl}
to obtain \ESlise{} filters and the generalized
\textsc{Gauss-Legendre} filters.  For \textsc{Zolotarev} filters, we
use the \href{http://guettel.com/rktoolbox/}{\tool{RKToolbox}} by
Güttel et al.  We will not consider other prominent examples of
rational filters, notably Trapezoid or former \SLISE{} filters, as
they are not as competitive.
Trapezoid filters offer a strictly monotonous decay in function value
for $\abbs{x} > 1$, which is not as sharp as for
\textsc{Gauss-Legendre} and \textsc{Zolotarev} filters leading to
significantly larger \WCR{} values~\cite{guettel}.  As for \SLISE{}
filters, Figure~\ref{fig:its_1.1} provides clear evidence that
\ESlise{} filters are superior when it comes to number of iterations
to convergence.
\subsection{Comparison of \WCR{} values}

\pgfplotstableread[col sep=comma]{figures/convergence_rates.dat}\data
\pgfplotstableread[col
sep=comma]{figures/convergence_rates_reduced.dat}\datar

\newcommand{\findmin}[4]{
	\pgfplotstablevertcat{\datatable}{#1}
	\pgfplotstablecreatecol[
	create col/expr={%
		\pgfplotstablerow
	}]{rownumber}\datatable
	\pgfplotstablesort[sort key={#2},sort cmp={float
	<}]{\sorted}{\datatable}%
	\pgfplotstablegetelem{0}{rownumber}\of{\sorted}%
	\pgfmathtruncatemacro#3{\pgfplotsretval+#4} 
	\pgfplotstableclear{\datatable}
}

\pgfplotstableset{
	highlight row min/.code 2 args={
		\pgfmathtruncatemacro\rowindex{#2-1}
		\pgfplotstabletranspose{\transposed}{#1}
		\findmin{\transposed}{\rowindex}{\maxval}{1}
		\edef\setstyles{\noexpand\pgfplotstableset{
				every row \rowindex\space column
				\maxval\noexpand/.style={
					postproc cell content/.append style={
						/pgfplots/table/@cell
						content/.add={$\noexpand\bf}{$}
					},
				}
			}
		}\setstyles
	},
}

\begin{table}
	\centering
	\tiny
	\pgfplotstabletypeset[
	columns={G,q,gauss,zolo,best}, 
	columns/G/.style={
		precision=6,
		column name={Gap $G$}
	},
	columns/q/.style={
		int detect,
		column name={Poles},
		preproc/expr={4*##1},
		column type/.add={}{|}
	},
	columns/q/.style={
		int detect,
		column name={\makecell{Poles per\\Quadrant $m$}},
		column type/.add={}{|}
	},
	columns/{gaussinf}/.style={
		sci,sci zerofill,sci e,precision=2,
		column name=\makecell{\textsc{Gauss-L.}\\$S=\infty$}
	},
	columns/{gauss}/.style={
		sci,sci zerofill,sci e,precision=2,
		column name=\makecell{\textsc{Gauss}-\\\textsc{Legendre} 
		},
	},
	columns/{zolo}/.style={
		sci,sci zerofill,sci e,precision=2,
		column name={\textsc{Zolotarev}}
	},
	columns/{best}/.style={
		sci,sci zerofill,sci e,precision=2,
		column name=\ESlise{}
	},
	every nth row ={5}{before row=\midrule},
	every head row/.style={before row=\toprule,after row=\midrule},
	every last row/.style={after row=\bottomrule},
	highlight row min ={\datar}{1},
	highlight row min ={\datar}{2},
	highlight row min ={\datar}{3},
	highlight row min ={\datar}{4},
	highlight row min ={\datar}{5},
	highlight row min ={\datar}{6},
	highlight row min ={\datar}{7},
	highlight row min ={\datar}{8},
	highlight row min ={\datar}{9},
	highlight row min ={\datar}{10},
	highlight row min ={\datar}{11},
	highlight row min ={\datar}{12},
	highlight row min ={\datar}{13},
	highlight row min ={\datar}{14},
	highlight row min ={\datar}{15},
	highlight row min ={\datar}{16},
	highlight row min ={\datar}{17},
	highlight row min ={\datar}{18},
	highlight row min ={\datar}{19},
	highlight row min ={\datar}{20},
	]{\data}

	\caption[Convergence rates]{\WCR{} values for different
	filters, gap parameters, and numbers of poles (smaller is better; 
	row-wise minimum in bold).
}
 \label{tab:convergence_rates}
\end{table}

In Table~\ref{tab:convergence_rates}, we list \WCR{} values for
\textsc{Gauss-Legendre}, \textsc{Zolotarev}, and \ESlise{} filters for
different gap parameters $G$ and poles per quadrant
$m$. As $m$ increases, \textsc{Zolotarev} filters feature a reliable, gradual decrease
in \WCR{}; by construction, they do not offer a decay
in function value as $\abbs{x} \rightarrow \infty$, unlike
\textsc{Gauss-Legendre} and \ESlise{}. Hence, \textsc{Gauss-Legendre}
and \ESlise{} filters lead to quicker convergence for some \FEAST{}
instances, even if their \WCR{} is larger. \textsc{Gauss-Legendre}
filters show the largest \WCR{}s and are not competitive with regards
to worst-case performance.  \ESlise{} exhibit a significant reduction
of \WCR{} compared to previous filters, especially for the default
choice of $m=4$ within the \FEAST{} eigensolver. Yet, the improvement
over \textsc{Zolotarev} filters diminishes as $m$ increases, because
the dimension of the underlying minimization increases with $m$. Large
numbers of poles per quadrant $m > 7$ correspond to high-degree
rational functions and are not taken into consideration. It is
important to notice that \ESlise{} filters seem to perform
particularly well for large $G$ and $m$, making them especially
competitive in the presence of eigenvalue clusters near the
interval boundaries.

The numbers in bold in Table~\ref{tab:convergence_rates} show that,
for almost all pairs $(G,m)$, \ESlise{} filters have the lowest \WCR{}
value. Based on Theorem~\ref{th:wcr}, the best worst-case performance
is offered by the filter with smallest \WCR{}. As we will see in the
following sections, this claim is confirmed by our numerical results.

\subsection{Experiments with \FEAST{}}
\label{sec:feastexperiments}
For our numerical experiments, we used \FEAST{} in the
version \tool{3.0},
compiled with the \tool{Intel Compiler 17.0.0},
and run on an \tool{Intel Core i7-6900K}. We selected the default \FEAST{}
parameters, with the exception of disabling\footnote{This behavior
  can be achieved through the \FEAST{} parameter
  $\mathtt{fpm(10)=1}$. For details, consult the \FEAST{}
  documentation \cite{feastsolver}.}  the repeated factorization of
the underlying linear systems. This substantially reduces runtime,
but requires sufficient \tool{RAM}.

A required argument of \FEAST{} is an upper bound $M_0$ on the number
of eigenvalues $M$ in the search interval; $M_0$ indicates the
size of the reduced eigenvalue problem, that is solved in every
\FEAST{} iteration.  As already mentioned in Sec.~\ref{sec:method},
the convergence rate of \FEAST{} substantially depends on this
subspace size $M_0$. Let us denote with $I \supset [-1,1]$ the interval
centered around $0$ that contains $M_0 > M$ eigenvalues;
then, FEAST{}'s convergence rate
is proportional to
\begin{equation}
	\label{eqn:convrate}
	\frac{\abs{r(\lambda_{\rm out})}}{\abs{r(\lambda_{\rm in})}},
\end{equation}
where $\abs{r(\lambda_{\rm in})} = \min_{\lambda\in [-1,1]\cap
  \sigma(A)} \abs{r(\lambda)}$, and $\abs{r(\lambda_{\rm out})} = \max_{\lambda\notin I\cap
  \sigma(A)} \abs{r(\lambda)}$.
A smaller subspace size $M_0$ yields faster \FEAST{} iterations, but
decreases the convergence rate and thus increases the number of
iterations.  As a compromise, the original \FEAST{} publication
\cite{feast} suggested a subspace size of
$M_0=\lceil C\times M \rceil$, for $C=1.5$.  This provides reliable
convergence within \FEAST{}, but not necessarily fastest convergence
as we see in the following.  Because only estimates of the actual
eigencounts are available in advance, we assess different scenarios by
studying a number of eigencount multipliers $C$.

It is important to notice that in all our comparison we change only
the rational filter and maintain all other part of \FEAST{}
unchanged. In particular we let \FEAST{} use the same default linear
system solves to tackle \eqref{eq:linsys} and direct eigensolver in
the Rayleigh-Ritz step. Since the level of accuracy of \FEAST{} is for
all practical purposes determined by these two tasks~\cite{feast_conv,
  feast_conv_2}, we consider the accuracy of the determined eigenpairs
across distinct filters comparable. In other words, since FEAST
reaches convergence with all filters using the same procedures, then
the solution are considered to be equally accurate.


\subsubsection{Experiment I}
\label{sec:exp-1}
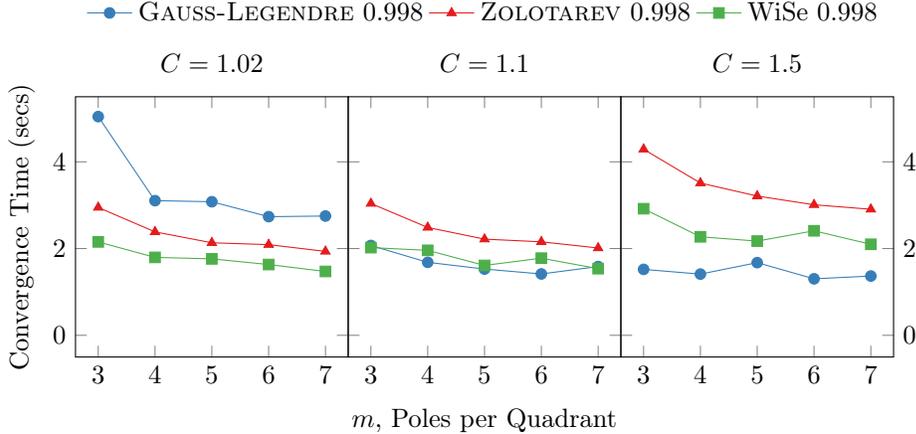
\begin{figure}
 \centering
 \tikzset{external/export next=false}
 \pgfplotstableread[col sep = comma]{figures/cnt_secs.csv}\cnt
 \begin{tikzpicture}
  \begin{groupplot}[group style = {
      horizontal sep = 0pt, vertical sep = 0pt, group size = 3 by 1},
    legend style={draw=none},
    ymax=5,
    ymin=0,
    enlargelimits=true,
    width = 0.4\textwidth, height = 0.24\textheight,
    xtick={3,...,7},
   ]
   \nextgroupplot [ title=\text{$C=1.02$}, 
   ylabel={Convergence Time (secs)}
]
	\addplot[color=GAUSSC,mark=*] table[x index={0},y index={1}] {\cnt};
   \addplot[color=ZOLOC,mark=triangle*] table[x index={0},y index={2}] 
   {\cnt};
   \addplot[color=SLISEC,mark=square*] table[x index={0},y index={3}] 
   {\cnt};

   \nextgroupplot [ title=\text{$C=1.1$}, 
   yticklabels={,,}]
   \addplot[color=GAUSSC,mark=*] table[x index={0},y index={4}] {\cnt};
   \addplot[color=ZOLOC,mark=triangle*] table[x index={0},y index={5}] 
   {\cnt};
   \addplot[color=SLISEC,mark=square*] table[x index={0},y index={6}] 
   {\cnt};
   
      \nextgroupplot [ title=\text{$C=1.5$}, 
      ylabel near ticks, yticklabel pos=right,
      legend style = { legend columns = -1, legend to name = grouplegend2}]
      
      \addplot[color=GAUSSC,mark=*] table[x index={0},y index={7}] {\cnt}; 
      \addlegendentry{\textsc{Gauss-Legendre} $0.998$}
   \addplot[color=ZOLOC,mark=triangle*] table[x index={0},y index={8}] 
   {\cnt}; \addlegendentry{\textsc{Zolotarev} $0.998$}
   \addplot[color=SLISEC,mark=square*] table[x index={0},y index={9}] 
   {\cnt}; \addlegendentry{\ESlise{} $0.998$}
  \end{groupplot}
  \node at ($(group c1r1.north)!0.5!(group c3r1.north)$) [above, 
  yshift=20] {\hypersetup{linkcolor=black}\ref{grouplegend2}};
    \node at ($(group c1r1.south)!0.5!(group c3r1.south)$) [below, 
  yshift=-1\pgfkeysvalueof{/pgfplots/every axis title shift}] 
  {\text{$m$}, Poles per Quadrant};
 \end{tikzpicture}
 \caption{Time required to solve the \tool{CNT} eigenproblem
   \cite{feast,guettel} through \FEAST{}, for different eigencount
   multipliers $C$, numbers of poles, and filters.  Averages of $10$
   executions are reported.  All the three filters lead to eigenpairs
   with the same level of accuracy.
}
 \label{fig:cnt_its_residuals}
\end{figure}


In Figure~\ref{fig:cnt_its_residuals}, we compare the convergence of
\FEAST{} for different eigencount multipliers $C>1$ on a specific
interval of the so-called \tool{CNT} eigenproblem. The matrices,
corresponding to this generalized eigenproblem, represent the
discretized Hamiltonian and Overlap operators of a physical system
studied in the context of a specific Density Functional Theory method
and have been used to analyze the worst-case performance of \FEAST{}
in previous publications \cite{feast,guettel}. These sparse \tool{CNT}
matrices $A,B \in \R^{12450\times12450}$, with $\num{86808}$ non-zero
entries, define the interior eigenproblem $Ax=\lambda B x$, where one
is interested in obtaining the $M=100$ eigenvalues in the interval
$[-65.0,4.96]$. The figure is divided into three quadrants, each
corresponding to one of three increasing eigencount multipliers
$C=1.02, 1.1, 1.5$. On the x-axis of each quadrant, the convergence
time (in seconds) of the \FEAST{} eigensolver, equipped with three distinct
filters, is plotted against increasing numbers of poles per quadrant
$m$ for each of these filters.
The value of the parameter $G$ is kept
fixed for all filters across the entire figure.

The result of this experiment demonstrates that best worst-case
convergence of the \FEAST{} eigensolver correlates strongly with the
\WCR{} value of the filter that is used to project onto the active
subspace. On the other hand, the size of the active subspace $M_0$
also contributes to the convergence time and can be a confounding factor in
interpreting the numerical results. In order to minimize the influence
of the latter on our interpretation of \WCR{}, we consider first
the quadrant for $C=1.02$. From Table~\ref{tab:convergence_rates}, we
expect that best convergence time is achieved by \FEAST{} when equipped with \ESlise{} filters, followed by \textsc{Zolotarev} and \textsc{Gauss-Legendre} filters, respectively. This is indeed the
case; the performance of the three filters is clearly separated by a
gap, which reflects the differences in \WCR{} between the filters for all $m$.

The influence of $m$ on filter performance
is less pronounced. This is not surprising. For instance, the
\textsc{Zolotarev} filter enables \FEAST{} to converge with a very slow
decrease of convergence rate as the number of poles increases, for
all considered $C$. The only consequence of increasing $C$ is a growth in
convergence time, since the size of the subspace increases with $C$ and
more operations with vectors must be performed by the linear system solver. This
behavior is due to the so-called equi-oscillation of
\textsc{Zolotarev} filters: These filters are optimal in approximating
the ideal filter in $\infty$-norm, but do not 
decay away from the filtered interval. \ESlise{} filters behave
similarly: While they decay (moderately) away from the interval, their
effectiveness is not determined by their value away from it but rather their
behavior very close to its boundary (see
Figure~\ref{fig:filters_log}). This is reflected by the very slight
decrease in time-to-convergence as $m$ increases although the
corresponding \WCR{} value decreases substantially as shown in Table~\ref{tab:convergence_rates}.

As the size of the active subspace increases, its influence on
convergence time becomes ever more pronounced. This is because the
difference between the convergence rate of \FEAST{} in
\eqref{eqn:convrate} may become increasingly larger than the \WCR{} of
the used filter as $C\gg 1.02$ \cite{guettel}. In other words, larger
active subspaces dilute the correlation between worst-case convergence
rate of the filter and convergence rate of the eigensolver. This is clearly
visible, if one traverses the quadrants in
Figure~\ref{fig:cnt_its_residuals} from left to right, and is best
illustrated by the \textsc{Gauss-Legendre} filters. These suffer from
slow convergence for small $C$, but can compensate for their large
\WCR{} for large subspace sizes because they decay rapidly in function
value away from the search interval. In other words, for
\textsc{Gauss-Legendre}, the size of the active subspace is much more
relevant than the \WCR{} of the filter.



This simple analysis, based on a very specific interior eigenvalue
problem, seems to suggest that \ESlise{} filters should always be
preferred over \textsc{Zolotarev} filters in those cases in which the
\ESlise{} \WCR{} is smaller than \textsc{Zolotarev} (compare Table~\ref{tab:convergence_rates}).
When comparing the performance of \ESlise{} and
\textsc{Gauss-Legendre} filters, it seems that the size of the active
subspace plays a major role in identifying the point 
at which one filter outperforms the other. In order to address this
question, we examine the convergence rate of \FEAST{} in the next
section, both in terms of number of FLOPs and number of subspace
iterations, using a large set of representative eigenproblems.

\subsubsection{Experiment II}
\label{sec:experiments_benchmarkset}

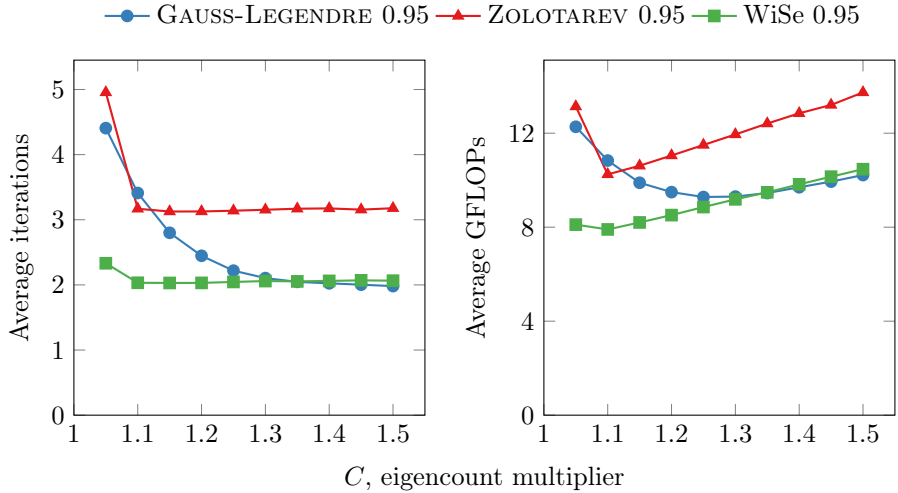
\begin{figure}
	\centering
	\tikzset{external/export next=false}
	\pgfplotstableread[col sep = comma]{figures/benchmark_set_0.95.dat}\bench
	\begin{tikzpicture}
		\begin{groupplot}[group style = {horizontal sep = 45pt, group size 
		= 2 by 1}, width = 0.48\textwidth, height = 0.3\textheight,legend 
		style={draw=none}]
		\nextgroupplot[ ymin=0, xmin=1.0,
						xmax=1.55,
						ytick distance=1,
						xtick={1.0,1.1,...,1.5,1.5},
						ylabel={Average iterations}]
			\addplot[color=GAUSSC,mark=*, thick] table[x index={0},y 
			expr=(\thisrowno{3}/2117)] {\bench}; 
			\addplot[color=ZOLOC,mark=triangle*, thick] table[x 
			index={0},y expr=(\thisrowno{5}/2117)] {\bench};
			\addplot[color=SLISEC,mark=square*, thick] table[x index={0},y 
			expr=(\thisrowno{1}/2117)] {\bench};
		\nextgroupplot [ymin=0, xmin=1.0,
						ytick distance=4,
						xmax=1.55,
						xtick={1.0,1.1,...,1.5,1.5},
						ylabel={Average GFLOPs},
						legend style = { legend columns = -1, legend to 
							name = grouplegend3}
						]
			\addplot[color=GAUSSC,mark=*, thick] table[x index={0},y 
			expr=(\thisrowno{4})] {\bench};
			\addlegendentry{\textsc{Gauss-Legendre $0.95$}}
			\addplot[color=ZOLOC,mark=triangle*, thick] table[x 
			index={0},y expr=(\thisrowno{6})] {\bench}; 
			\addlegendentry{\textsc{Zolotarev $0.95$}}
			\addplot[color=SLISEC,mark=square*, thick] table[x index={0},y 
			expr=(\thisrowno{2})] {\bench}; \addlegendentry{\ESlise{} 
			$0.95$}
		\end{groupplot}
	\node at ($(group c1r1.north)!0.5!(group c2r1.north)$) [above, 
		yshift=5] {\hypersetup{linkcolor=black}\ref{grouplegend3}};
	\node at ($(group c1r1.south)!0.5!(group c2r1.south)$) [below, 
	yshift=-1\pgfkeysvalueof{/pgfplots/every axis title shift}] {$C$, 
	eigencount multiplier}; 
	\end{tikzpicture}
	\caption{The average numbers of iterations and FLOPs (floating-point 
	operations),
	required by \FEAST{} to solve $2117$ benchmark problems \cite{jan}, 
	for different subspace sizes
	(by multiplying the actual eigencounts with a fixed scalar $C$)
	and filters with $G=0.95$.
	}
	\label{fig:benchmark_set_0.95}
\end{figure}

We consider a moderate gap parameter of $G=0.95$ and
a set of $\num{2117}$ interior eigenproblems.
These eigenproblems were obtained from \tool{Si2}, a sparse and symmetric
matrix from the University of Florida Matrix Collection
\cite{Si2}, by selecting $2117$ different search intervals
$[a,b]$\footnote{The code to obtain such benchmark sets from arbitrary
  matrices is freely available at
  \url{https://github.com/SimLabQuantumMaterials/SpectrumSlicingTestSuite.jl}.}
as described in \cite[Appendix B]{jan}. Each search interval uniquely
identifies an interior eigenproblem with its unique eigenvalue
distribution and eigenvalue count. As such, it is quite general
and statistically relevant, since it reflects the large variations
that are possible in distributing and clustering eigenvalues inside,
outside, and in the vicinity of the search interval ends.

We initially solved for each of the $\num{2117}$ benchmark problems
with a fixed value of $m=4$ poles per quadrant, the default
in the \FEAST{} eigensolver. As in the previous section, we repeated this test for all
three filters for increasing values of
the eigencount multiplier $C$. The results of these tests are
graphically reported in Figure~\ref{fig:benchmark_set_0.95}, which
plots the number of subspace iterations and the total number of
FLOPs performed by \FEAST{}. These results confirm the analysis of
Section~\ref{sec:exp-1}.  The \textsc{Zolotarev} and \ESlise{} filters
maintain a linear behavior as a function of increasing dimension of
the active subspace as soon as $C\ge1.1$. In other words, the \WCR{}
of these filters influences the convergence of the eigensolver only for
active subspaces that closely match the true number of eigenvalues in
the interval $[a,b]$. As soon as the size of the active subspace gets
larger, the convergence of the eigensolver is dictated by
\eqref{eqn:convrate}. This interpretation is made even clearer when
one looks at the linear increase in FLOPs for \FEAST{}, equipped with
these two filters: While the average number of iterations remains
constant, the total number of floating-point operations increases due
to the linear increase in the total number of right-hand-side vectors
$Y$ for which \eqref{eq:linsys} must be solved.

As $C$ grows, the rate of convergence of \FEAST{}, equipped with the
\textsc{Gauss-Legendre} filter, equals the one of \ESlise{}---which is
dictated by the \WCR{} at smaller values of $C$.  For subspace
multipliers $C > 1.3$, the \textsc{Gauss-Legendre} starts competing,
on average, with the \ESlise{}. It must be noted that this comes at a
cost: Using \textsc{Gauss-Legendre} for a relatively large active
subspace, such as the default value of $C=1.5$ suggested by \FEAST{},
has on average a higher FLOP count than the \ESlise{} filters for
eigenvalue counts $C<1.3$.
This observation is fairly independent from the number of
poles per quadrant used, as shown in
Figure~\ref{fig:benchmark_set_morepoles}.  \textsc{Gauss-Legendre}
filters have a slight advantage with respect to FLOP count for large
subspace sizes ($C=1.5$), but they behave worse than the \ESlise{} for
any number of poles and small subspace sizes ($C=1.1$). Due to the
decay in value of the filter function, the \ESlise{} filter even
outperforms \textsc{Zolotarev} for $m=7$, despite its larger
\WCR{}. For larger $m$, overall FLOP count increases, while the
differences in FLOP count across the filters shrinks.

In conclusion, \FEAST{}, equipped with \ESlise{} filters, offers a
competitive advantage over the use of \textsc{Gauss-Legendre} and
\textsc{Zolotarev} filters. \ESlise{} filters are quite stable with
respect to the convergence rate of the eigensolver, irrespective of the
active subspace or the degree of the filter function. Their use seems
to almost always minimize the total FLOP count required by \FEAST{} to
reach convergence. In addition, their effectiveness for small
eigencount multipliers suggests that \ESlise{} filters should be
preferred in all those cases where it is necessary to contain the
subspace size, either because the \tool{RAM} is limited or the
underlying spectrum distribution is unknown.

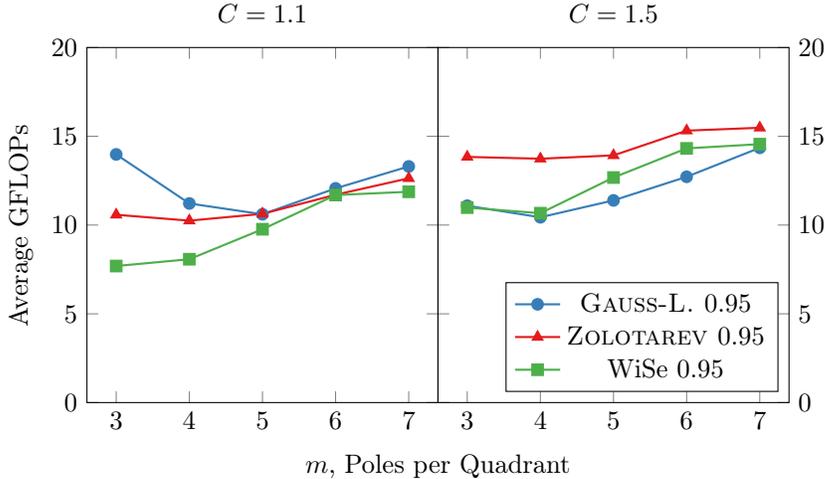
\begin{figure}
	\centering
	\tikzset{external/export next=false}
	\pgfplotstableread[col sep = 
	comma]{figures/benchmark_set_morepoles.csv}\bench
	\begin{tikzpicture}
		\begin{groupplot}[group style = {horizontal sep = 0pt, group size 
		= 2 by 1}, width = 0.48\textwidth, height = 0.3\textheight]
		\nextgroupplot[ ymin=0,
						ymax=20,
						title=\text{$C=1.1$},
						xtick={3,...,7},
						ylabel={Average GFLOPs},
						legend style = { legend columns = -1, legend to 
						name = grouplegend}]
			\addplot[color=GAUSSC,mark=*, thick] table[x index={0},y 
			index={1}] {\bench}; 
			\addplot[color=ZOLOC,mark=triangle*, thick] table[x 
			index={0},y index={2}] {\bench}; 
			\addplot[color=SLISEC,mark=square*, thick] table[x index={0},y 
			index={3}] {\bench}; 
		\nextgroupplot [ymin=0,
						ymax=20,
						title=\text{$C=1.5$},
						xtick={3,...,7},
						legend pos=south east,
						ylabel near ticks, yticklabel pos=right,
						]
			\addplot[color=GAUSSC,mark=*, thick] table[x index={0},y 
			index={4}] {\bench};
			\addplot[color=ZOLOC,mark=triangle*, thick] table[x 
			index={0},y index={5}] {\bench};
			\addplot[color=SLISEC,mark=square*, thick] table[x index={0},y 
			index={6}] {\bench};
			\legend{\textsc{Gauss-L. $0.95$},\textsc{Zolotarev 
			$0.95$},\ESlise{} $0.95$}
		\end{groupplot}
	\node at ($(group c1r1.south)!0.5!(group c2r1.south)$) [below, 
	yshift=-1\pgfkeysvalueof{/pgfplots/every axis title shift}] {$m$, 
	Poles per Quadrant}; 
	\end{tikzpicture}
	\caption{Comparison of average FLOPs, required to solve $2117$ 
	benchmark problems \cite{jan}, using \textsc{Gauss-Legendre}, 
	\textsc{Zolotarev}, and \ESlise{} filters for different of poles 
	numbers $m$ per quadrant,
	while $C \in \set{1.1, 1.5}$ and $G=0.95$.
	}
	\label{fig:benchmark_set_morepoles}
\end{figure}


\section{Conclusions}

In this work, we show how we decreased time to convergence of the
\SLISE{} optimization framework by using in it the minimization
algorithm \LBFGSB{}.  When computing a box-constrained \SLISE{}
filter, only hundreds of function evaluations are needed instead of
millions. We exploit the improved performance by introducing a second
optimization process for the numerical minimization of the Worst-case
Convergence Rate of the \SLISE{} rational filters. The byproduct of
such minimization is the elimination of the dependence of the filter
from the weight functions used in the non-linear least squares
functional. 
The new \ESlise{} filters outperform \textsc{Gauss-Legendre} and
\textsc{Zolotarev} filters both in terms of execution time, number of
subspace iterations, and FLOPs count necessary to reach
convergence by the eigensolver.



Increasing the performance of the optimization of rational filters and
eliminating their dependence from a number of adjustable parameters
has an additional indirect and important impact on the eigensolver
using the rational filters as spectral projectors. This class of
solvers lend themselves to multiple levels of parallelism: At the
highest level each interval $[a,b]$ can be split in subintervals
$[a_j,b_j]$, each of which constitutes a trivially separate eigenproblem; at a
mid-level the spectral solver requires the solution of a linear system
for each pole $z_i$; at the lowest level each linear system has to be
solved for multiple RHS. While such a general scheme makes this class
of eigensolver attractive, it complicates substantially the problem of
balancing the computational load. One of the main contributor to the
uncertainty of a well-balanced computation is the ability of the
spectral filter in determining the number of subspace iterations
needed to converge the full subspace corresponding to each
$[a_j,b_j]$. 

Our \ESlise{} filters overcome this uncertainty by: 1) decoupling the rate
of convergence from the size of the active search subspace, 2) drastically
reducing the dependence of the number of poles which can be safely set
to a standard value (e.g. $m=4$ in FEAST).
The net result is that the spectral filter has the same effectiveness
for any sub-interval selected: For a given linear system solver the
number of iterations required to reach convergence is minimized and
independent from the eigenvalues distribution. Load balancing is then
achieved by choosing sub-intervals with approximately the same
eigenvalue count. Since obtaining a good estimate for the eigenvalue
count and the eigenvalue distribution is a solved
problem~\cite{eigencount, linlin}, the result presented
in this paper eliminates the influence of the spectral filter on load
balancing for all practical purposes. The remaining challenge is
balancing the load when solving for distinct linear system with
multiple RHS. This is the focus of further ongoing work.





\section*{Acknowledgements}
We thank Jan Winkelmann for having provided support to the first
author in developing the bulk of the work that contributed to this
paper and Sebastian Achilles for useful discussions.

\bibliographystyle{siamplain}
\bibliography{bibliography}
\end{document}